\documentclass[basic]{custom}
 \usepackage{mathptmx, xcolor}
 \usepackage{latexsym, amsmath, amssymb, bm, mathtools} 
 \usepackage{hyperref, multirow} \usepackage[authoryear]{natbib}
 \usepackage{float}

 \def\NI{\noindent}
 \def\TS{\textstyle} \def\d{{\rm d}} \def\PHI{\varphi}
 \def\VEC#1{{\pmb{#1}}} \def\MAT#1{{\pmb{#1}}}
 \def\IE{{\it i.e.}} \def\EG{{\it e.g.}}  
 \def\FIRST{{$1^{\textrm{st}}$}} \def\SECOND{{$2^{\textrm{nd}}$}} 
 \def\THIRD{{$3^{\textrm{rd}}$}} \def\NTH#1{{${#1}^{\textrm{th}}$}} 
 \def\DOTS{{...}} \def\DS{\displaystyle}

 \def\STRUT{\vphantom{\vert_\vert^\vert}}

 \def\TRUE{$\textsf{T}$} \def\FALSE{$\tiny\textsf{F}$}
 \def\DOI#1{\href{DOI: #1}{https://doi.org/#1}}

 \def\LEAF{\includegraphics[height=19pt]{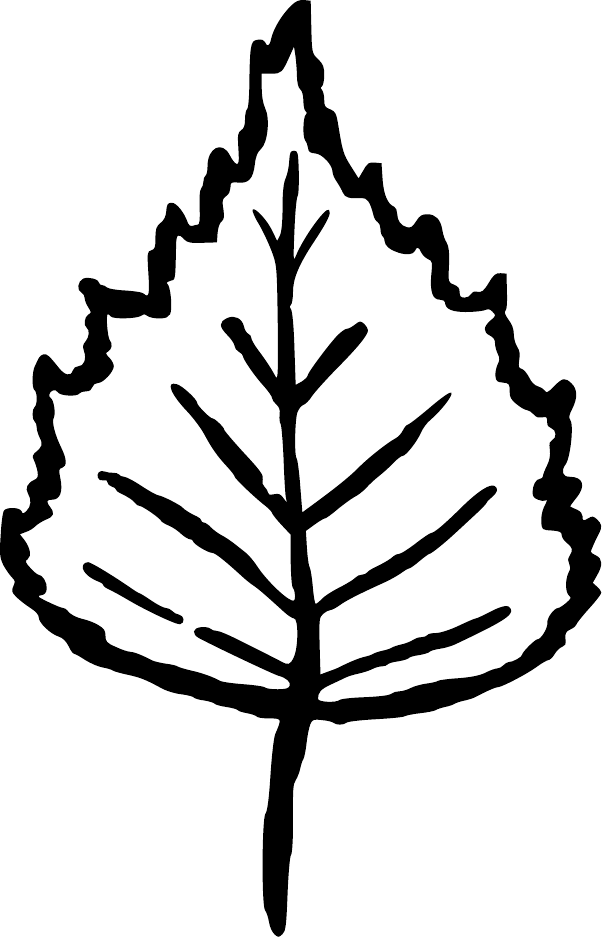}}
 \def\ROOT{\includegraphics[height=22pt]{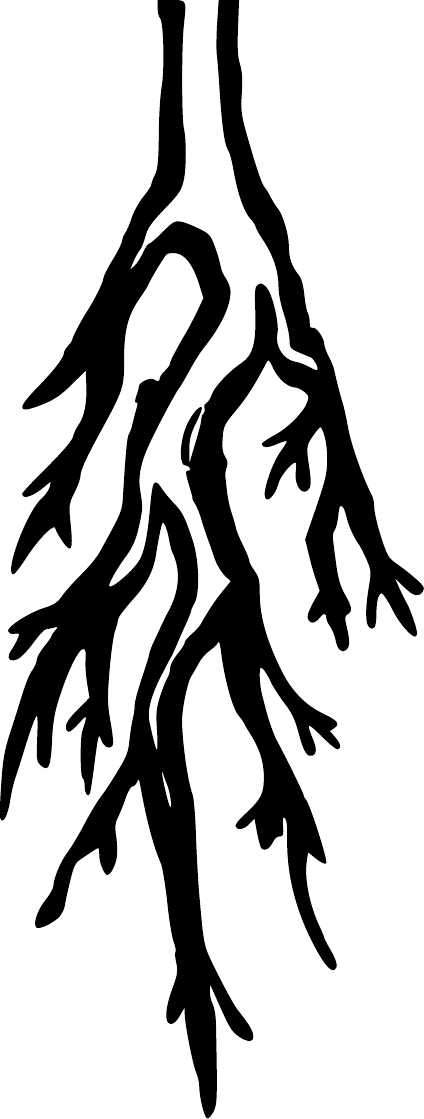}}
 \def\BIRCH{\includegraphics[height=49pt]{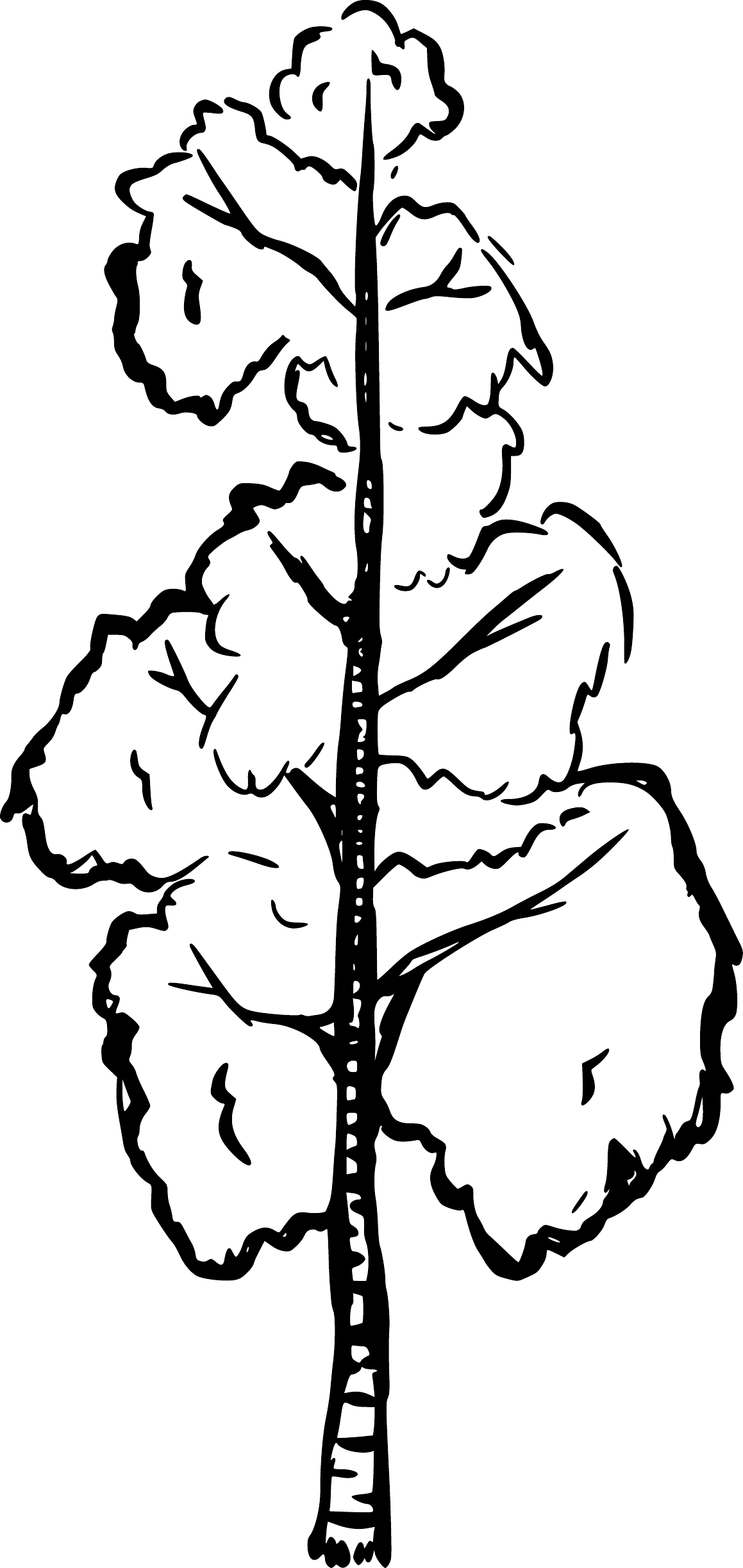}}
 \def\LARCH{\includegraphics[height=49pt]{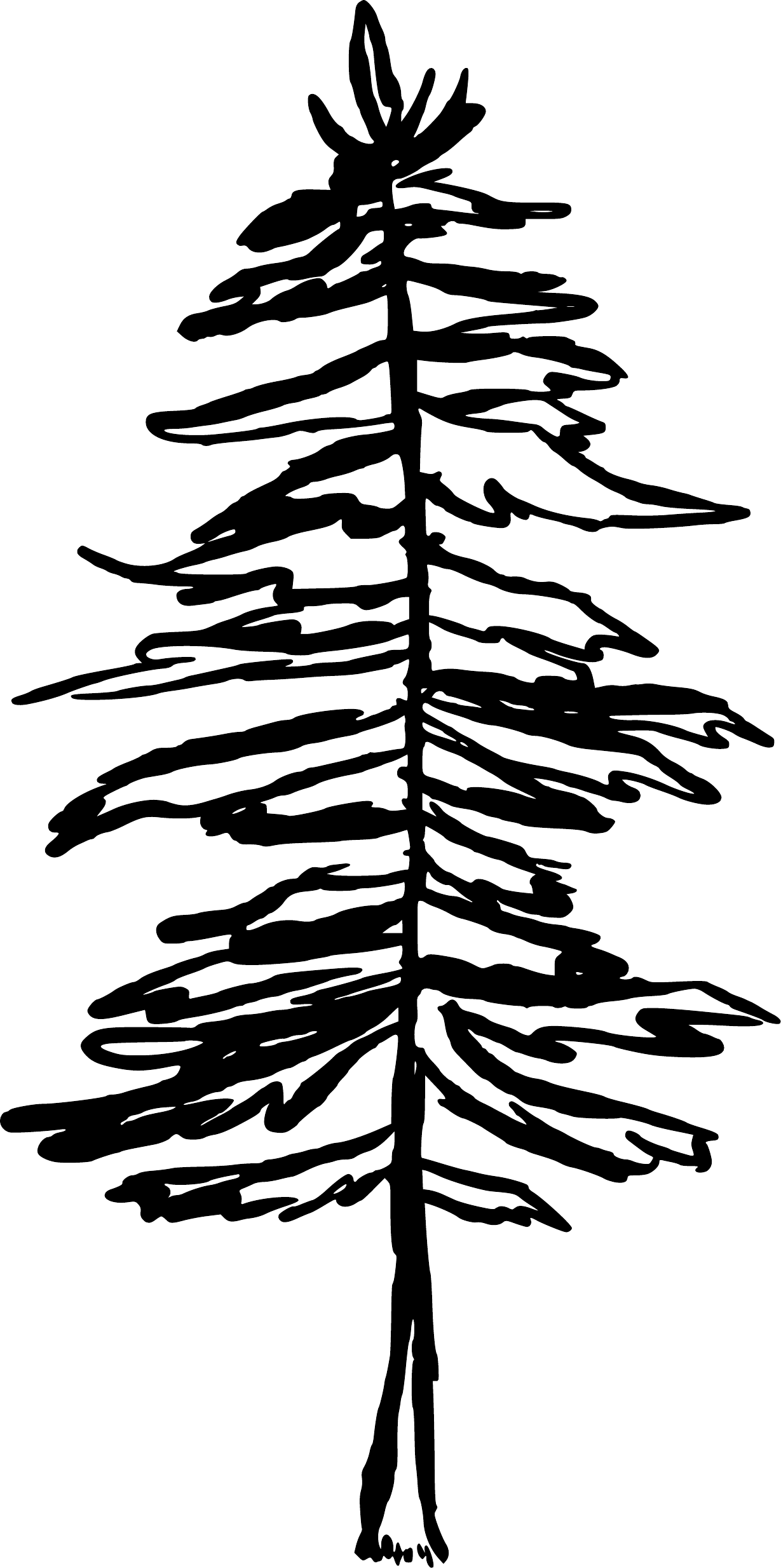}}

 \def\SEP{,\mkern2mu} \def\DWS{\mkern3.5mu}
 \def\AIJ{{a_{\mkern0.85mu i \mkern-0.85mu j}}}
 \def\BI{{b_{\mkern0.8mu i}}}
 \def\BJ{{b_{\mkern-1.25mu j}}}
 \def\CI{{c_{\mkern0.5mu i}}}
 \def\FI{{\VEC{F}_{\mkern-1.5mu i}}}
 \def\FJ{{\VEC{F}_{\mkern-4mu j}}}
 \def\MIJ{{m_{\mkern0.85mu i \mkern-0.85mu j}}}
 \def\DT{\d \mkern1mu t}
 \def\X12S#1{{\VEC{X}_{\mkern-1mu{#1}}}}
 \def\APHI{{\MAT{A} \mkern0.5mu \VEC{\Phi}}}
 \def\BPHI{{\VEC{b} \mkern2mu \VEC{\Phi}}}
 \def\BRAT{( \mkern0.25mu t \mkern-0.5mu )}
 \def\BRA1{(\VEC{1} \mkern-0.75mu )}
 \def\BRAA1{(\MAT{A} \VEC{1} \mkern-0.75mu )}
 \def\BRAP1{(\VEC{p}_1 \mkern-0.5mu )}
 \def\BRAAP1{(\MAT{A} \VEC{p}_1 \mkern-0.5mu )}
 \def\BRAAQ1{(\MAT{A} \mkern0.5mu \VEC{q}_1 \mkern-0.5mu )}
 \def\INT01{$[ \mkern1mu 0, 1 \mkern0.5mu ]$}
 \def\OMEGO{\omega_{\mkern1mu 1}}
 \def\OMEGT{\omega_{\mkern1.5mu 2}}
 \def\OMEGH{\omega_{\mkern1.5mu 3}}

\begin{document}

\title{Eight-stage pseudo-symplectic Runge--Kutta methods of order 
(4,\;8)}

\author*{\hskip42pt\fnm{Misha} \sur{Stepanov}
\qquad {\small stepanov@arizona.edu}}

\affil{\orgdiv{Department of Mathematics \,\;and\;\, Program in Applied 
Mathematics}, \orgname{University of Arizona}, 
\orgaddress{\city{Tucson}, \state{AZ 85721}, 
\country{US$\mkern0.5mu$A}}}

\abstract{\hskip10pt{}Using simplifying assumptions that are related to 
the time reversal symmetry, a 1-dimensional family of $8$-stage 
pseudo-symplectic Runge--Kutta methods of order $(4 \SEP 8)$, \IE, 
methods of order $4$ that preserve symplectic structure up to order $8$, 
is derived. An example of $7$-stage method of order $(4 \SEP 9)$ is 
given.}

\keywords{pseudo-symplectic Runge--Kutta methods}

\pacs[MSC Classification]{65L05, 65L06}

\maketitle

\vspace{-12pt}

\hskip\parindent{}Runge--Kutta methods (see, \EG, \citep[s.~23 and 
ch.~3]{But16}, \citep[ch.~II]{HNW93}, \citep[ch.~4]{AsPe98}, 
\citep[ch.~3]{Ise08}) are widely and successfully used to solve ordinary 
differential equations numerically for over a century \citep[]{BuWa96}. 
Being applied to a system $\d \mkern0.5mu \VEC{x} / \DT = \VEC{f}(t \SEP 
\VEC{x})$, in order to propagate by the step size $h$ and update the 
position, $\VEC{x}\BRAT \mapsto \VEC{\tilde{x}}(t + h)$, where 
$\VEC{\tilde{x}}(t + h)$ is a numerical approximation to the exact 
solution $\VEC{x}(t + h)$, an \mbox{$s$-stage} Runge--Kutta method 
(which is determined by the coefficients $\AIJ$, weights $\BJ$, and 
nodes $c_{\mkern0.5mu i}$) would form the following system of equations 
for $\X12S{1}$, $\X12S{2}$, \DOTS, $\X12S{s}$:

\vspace{-11pt}

 \begin{gather*}
   \X12S{i} = \VEC{x}\BRAT + h \mkern-1.5mu \sum_{j = 1}^{s} \AIJ 
\mkern2mu \FJ , \qquad \FI = \VEC{f} \bigl( t + \CI \mkern1mu h \SEP 
\X12S{i} \bigr), \qquad i = 1, \mkern1.5mu 2, \DOTS, \mkern1.5mu s
 \end{gather*}

\vspace{-5pt}

\NI{}solve it, and then compute $\VEC{\tilde{x}}(t + h) = \VEC{x}\BRAT + 
h \sum_{\mkern-1mu j = 1}^s \BJ \mkern2mu \FJ$. In the limit $h \to 0$ 
all the vectors $\FI$, where $1 \le i \le s$, are the same, so it is 
natural and will be assumed that $\sum_{\mkern-1mu j = 1}^{s} \AIJ = 
\CI$ for all $i$.

It is a common wisdom that whenever a numerically simulated system has a 
conservation law, in order to get sensible results it is often desirable 
to use so-called conservative schemes, \IE, the ones that preserve the 
conserved quantity exactly (up to round-off errors). A practically 
important case is when the system is Hamiltonian \citep[ch.~8]{Arn89}, 
\citep[s.~VI.1]{HLW06}, \citep[p.~251]{But21}. Integration methods are 
called symplectic if they conserve the symplectic structure.

Let $\MAT{M} = \bigl[ \MIJ \bigr]$ be the $s \times s$ matrix with $\MIJ 
= \BI \mkern1.2mu \AIJ + \BJ \mkern1.5mu a_{\mkern-1.2mu j \mkern1.2mu 
i} - \BI \mkern1.2mu \BJ$ as its matrix element in the $i^{\mkern1mu 
\textrm{th}}$ row and $j^{\mkern1mu \textrm{th}}$ column (see, \EG, 
\citep[]{But75}, \citep[eq.~(2.2)]{BuBu79}, \citep[]{Coo87}, 
\citep[s.~390]{But16}, \citep[p.~316]{HNW93}, \citep[p.~79]{Ise08}, 
\citep[p.~253]{But21}). A Runge--Kutta method is symplectic if $\MAT{M}$ 
is a zero matrix $\MAT{O}$, and only if it is equivalent to a method 
with $\MAT{M} = \MAT{O}$, as was independently shown by several authors: 
\citep[]{San88}, \citep[]{Sur88}, \citep[]{Sur89}, 
\citep[]{Las88}.\footnote{~An $S$-reducible (see, \EG, 
\citep[p.~188]{HW96}) $2$-stage method with $a_{11} = a_{22} = b_1 = b_2 
= c_1 = c_2 = \smash{\frac12}$ and $a_{12} = a_{21} = 0$ has $m_{11} = 
m_{22} = -m_{12} = -m_{21} = \smash{\frac14} \ne 0$. It is equivalent to 
the midpoint method and is symplectic.}

Unless a Hamiltonian is of a special type, \EG, as it often happens in 
problems coming from classical mechanics, is separable: 
$\mathcal{H}(\VEC{p} \SEP \VEC{x}) = T(\VEC{p}) + U(\VEC{x})$, where 
$\VEC{x}$ and $\VEC{p}$ are canonical coordinates, symplectic methods 
are implicit. In \citep[]{AuCh98a} the concept of so-called 
pseudo-symplectic methods was introduced. A method is said to be of 
pseudo-symplectic order $(p \SEP q)$ if it is of order $p$ (for order 
conditions see, \EG, \citep[p.~172]{But16}, \citep[p.~126]{But21}, 
\citep[p.~148]{HNW93}), and the symplectic structure is conserved up to 
the order $q$ (see, \EG, \citep[eqs.~(2.5), (2.7), and 
tab.~2.1]{AuCh98a}).\footnote{~In \citep[eq.~(2.7)]{AuCh98a} ``$(t \SEP 
t') \in \hat{T}(k) \times \hat{T}(k)$'' should be read as ``$(t \SEP t') 
\in \hat{T}(k)$''.} It was shown that methods with $q \ge 2 p$ have 
better Hamiltonian conservation properties \citep[thm.~2.6]{AuCh98a}, 
and an explicit $5$-stage pseudo-symplectic Runge--Kutta method of order 
$(3 \SEP 6)$ was constructed \citep[fig.~4.1]{AuCh98a}.

With an abundance of implicit symplectic Runge--Kutta methods, it is 
only explicit pseudo-symplectic methods that are of practical interest. 
There is no $6$-stage method of order $(4 \SEP 8)$ 
\citep[prop.~4.2]{AuCh98a} (see also \citep[thm.~3.3]{AuCh98b}), but 
there are known $6$-stage methods of order $(4 \SEP 7)$: 
\citep[p.~262]{CLMR10} and \citep[p.~90]{CCRL17}.

The aim of this work is to construct a pseudo-symplectic method of order 
$(4 \SEP 8)$. Order conditions are discussed in 
Section~\ref{order_conditions}. A $1$-dimensional family of $8$-stage 
pseudo-symplectic methods of order $(4 \SEP 8)$ is derived in 
Section~\ref{family}. The efficiency of new methods is tested on three 
numerical examples in Section~\ref{numerical_tests}.

\section{Order conditions} \label{order_conditions}

\hskip\parindent{}The element-wise product of column vectors $\VEC{x}$ 
and $\VEC{y}$ will be denoted as $\VEC{x} \mkern0.5mu . \mkern1mu 
\VEC{y}$, \IE, $(\VEC{x} \mkern0.5mu . \mkern1mu \VEC{y})_i = x_i 
\mkern1mu y_i$. The element-wise product of $n$ copies of vector 
$\VEC{x}$ will be written as $\VEC{x}^{\mkern1mu n}$. Let $\VEC{1}$ be 
the $s$-dimensional column vector with all components being equal to 
$1$; $\MAT{A} = \bigl[ \AIJ \bigr]$ be the $s \times s$ matrix with 
$\AIJ$ as its matrix element in the $i^{\mkern1mu \textrm{th}}$ row and 
$j^{\mkern1mu \textrm{th}}$ column; $\VEC{b} = \bigl[ \BJ \bigr]$ be the 
weights row vector; and $\VEC{c} = \bigl[ \CI \bigr]$ be the nodes 
vector. Let $\VEC{q}_n = \MAT{A} \mkern1mu \VEC{c}^{\mkern0.5mu n} - 
\frac{1}{n + 1} \mkern1mu \VEC{c}^{\mkern0.5mu n + 1}$, \EG, $\VEC{q}_0 
= \VEC{0}$ as $\MAT{A} \VEC{1} = \VEC{c}$.

Given rooted trees $\textrm{t}_1$, $\textrm{t}_2$, \DOTS, 
$\textrm{t}_n$, a new tree $[ \mkern1mu \textrm{t}_1 \mkern2mu 
\textrm{t}_2 \mkern2mu \DOTS \mkern2mu \textrm{t}_n \mkern1mu ]$ is 
obtained by connecting with $n$ edges their roots to a new vertex, the 
latter becomes a new root \citep[s.~301]{But16}, \citep[p.~152]{HNW93}, 
\citep[p.~44]{But21}, \citep[p.~53]{HLW06}. Consider a vector function 
$\VEC{\Phi} : \mathrm{T} \to \pmb{\mathbf{R}}^s$ on the set of rooted 
trees that is recursively defined as $\VEC{\Phi}( 
\begin{picture}(4,6)(0,0) \put(3.35,2.5){\makebox(0,0){\circle*{3}}} 
\end{picture} ) = \VEC{1}$ and $\VEC{\Phi} \bigl( [ \mkern1mu 
\textrm{t}_1 \mkern2mu \textrm{t}_2 \mkern2mu \DOTS \mkern2mu 
\textrm{t}_n \mkern1mu ] \bigr) = \prod_{\mkern1mu m = 1}^{\mkern1mu n} 
\APHI(\textrm{t}_m)$, where the product of vectors is taken 
element-wise. This function coincides with \emph{derivative weights} 
\citep[def.~312A]{But16}, \citep[pp.~148 and 151]{HNW93}, 
\citep[def.~2.2]{AuCh98a}, it is closely related to \emph{internal} or 
\emph{stage weights} $\APHI(\textrm{t})$ \citep[p.~125]{But21} and 
\emph{elementary weights} $\BPHI(\textrm{t})$ \citep[p.~55]{HLW06}.

Representation of analytical expressions by pictures or diagrams is 
often convenient and illustrative, and is used in many different 
contexts: countless diagram techniques in theoretical physics; group 
theory (see, \EG, \citep[]{Cvi08}); graphical representation of Fredholm 
series \citep[s.~11.3]{MaWa70}; and connection of derivatives with 
rooted trees \citep[]{Cay57}, which was heavily developed in the context 
of Runge--Kutta methods \citep[]{Mer57}, \citep[]{But72}, 
\citep[ch.~3]{But16}, \citep[s.~II.2]{HNW93}, \citep[ch.~2]{But21}, 
\citep[ch.~III]{HLW06}, with occasional usage of picturesque trees in 
diagrams \citep[pp.~188 and 189]{But16}, \citep[pp.~175 and 328]{HNW93}.

Similarly to Penrose graphical notation \citep[]{Pen71}, an object that 
is a tensor of type $(m \SEP n)$ will be depicted by a symbol of the 
object connected by lines to $m$ and $n$ black dots below and above it, 
respectively:$\mkern1mu$\footnote{~Such an 
up$\mkern2mu$/$\mkern0.75mu$down orientation is the reverse of the one in 
\citep[]{Pen71}, and is used here so that a tree has its root at the 
bottom. Not drawing the dots allows an elegant representation of 
Kronecker delta $\delta^{\mkern1mu a}_{\mkern1mu b}$ and metric tensors 
$g^{\mkern1mu ab}$ and $g_{\mkern1mu ab}$ \citep[pp.~226 and 231]{Pen71}, 
while here the dots correspond to the vertices of rooted 
trees.}${}^{\SEP}$\footnote{~The drawings of a leaf, a root, of trees 
\emph{Betula pendula} and \emph{Larix sibirica} are by Olga Stepanova, 
used with permission.}
 \begin{align*}
   \begin{tabular}{r|l}
   \begin{picture}(17,23)(0,0)
     \put(0.9,2.5){\LEAF}
     \put(8,1.5){\makebox(0,0){\circle*{3}}}
   \end{picture} & \raisebox{3.67pt}{vector $\VEC{1}$} \\
   \hline
   \begin{picture}(60,55)(0,0)
     \put(1.5,44){\makebox(0,0){$\textrm{t}$}}
     \put(2.55,1.95){\BIRCH}
     \put(15,1.5){\makebox(0,0){\circle*{3}}}
     \put(34,1.5){\LARCH}
     \put(48,1.5){\makebox(0,0){\circle*{3}}}
   \end{picture}
     & \raisebox{12pt}{$\VEC{\Phi}(\textrm{t})$ or an arbitrary vector} \\
   \hline
   \begin{picture}(21,24)(0,0)
     \put(1.5,10){\makebox(0,0){\rotatebox{90}{\rule[3pt]{17pt}{1.44pt}}}}
     \put(1.5,1.5){\makebox(0,0){\circle*{3}}}
     \put(1.5,18.5){\makebox(0,0){\circle*{3}}}
     \put(16.5,10){\makebox(0,0){\rotatebox{90}{\rule[3pt]{17pt}{1.44pt}}}}
     \put(16.5,10){\makebox(0,0){\rotatebox{90}{\color{white}\rule[3pt]{17pt}{0.72pt}}}}
     \put(16.5,1.5){\makebox(0,0){\circle*{3}}}
     \put(16.5,18.5){\makebox(0,0){\circle*{3}}}
   \end{picture}
     & \raisebox{3.67pt}{matrix $\MAT{A}$ and identity matrix $\MAT{I}$} \\
   \hline
   \begin{picture}(16,29)(0,0)
    \put(0.4,0.5){\ROOT}
     \put(5.85,23.5){\makebox(0,0){\circle*{3}}}
   \end{picture}
     & \raisebox{6.67pt}{weights row vector $\VEC{b}$} \\
   \hline
   \begin{picture}(22,17)(0,0)
     \put(11,1.5){\makebox(0,0){\circle*{3}}}
     \put(6,7){\makebox(0,0){\rotatebox{135}{\rule[3pt]{14pt}{0.72pt}}}}
     \put(15.5,5.5){\makebox(0,0){\rotatebox{45}{\rule[3pt]{14.5pt}{0.72pt}}}}
     \put(8.75,7.4){\makebox(0,0){\rotatebox{112.5}{\rule[3pt]{12pt}{0.72pt}}}}
     \put(11.75,10.5){\makebox(0,0){\footnotesize$...$}}
   \end{picture}
     & \raisebox{1.33pt}{element-wise product of vectors}
   \end{tabular} \end{align*} Lines coming both from below and from 
above to a dot produce a tensor contraction, \IE, the summation over the 
values of the corresponding tensor index. Outer products are obtained by 
simply drawing objects next to each other \citep[p.~224]{Pen71}.

In bra-ket notation \citep{Dir39} vectors, linear forms, scalar and 
outer products (the latter is an example of an operator) are denoted 
through $\vert u \rangle$, $\langle f \vert$, $\langle f \vert u 
\rangle$, and $\vert u \rangle \langle f \vert$, respectively, with 
clear distinction between vectors and forms. In quantum mechanics the 
evolution is \emph{linear}, so there is no need for branching. 
Similarly, the stability function $R(z) = 1 + z \mkern2mu \VEC{b} 
\mkern0.5mu (\MAT{I} - z \mkern2mu \MAT{A})^{-1} \VEC{1} = 1 + z 
\mkern2mu \VEC{b} \mkern0.5mu \VEC{1} + z^{\mkern1mu 2} \VEC{b} 
\mkern2mu \MAT{A} \VEC{1} + z^{\mkern1mu 3} \VEC{b} \mkern2mu \MAT{A}^2 
\VEC{1} + \DOTS = 1 + \sum_{n = 0}^\infty z^{\mkern1mu n + 1} \VEC{b} 
\mkern2mu \MAT{A}^n \VEC{1}$ (see, \EG, \citep[s.~238]{But16}, 
\citep[s.~5.3]{But21}, \citep[s.~4.4]{AsPe98}), which describes the 
behavior of a Runge--Kutta method applied to a \emph{linear} equation 
$\d \mkern0.5mu x / \DT = \lambda \mkern1mu x$ with $z = \lambda 
\mkern1mu h$ (\IE, $\partial^2 \mkern-1mu f_{\mkern1mu i} / \partial 
x_{\mkern-1mu j} \mkern1.5mu \partial x_k = 0$), contains only the terms 
that correspond to rooting trees with no branching.

For the system $\d \mkern0.5mu \VEC{x} / \DT = \VEC{f}(t \SEP \VEC{x})$ 
its integral form

\vspace{-16pt}

 \begin{gather*}
   \VEC{x}(t + h) = \, \VEC{x}\BRAT \, + \int\limits_t^{t + h} \DT' 
\DWS \VEC{f} \bigl( t' \! \SEP \VEC{x}(t') \bigr) = \, \VEC{x}\BRAT \, + 
\, h \int\limits_0^1 \d \theta \DWS \VEC{f} \bigl( t + \theta h \SEP 
\VEC{x}(t + \theta h) \bigr)
 \end{gather*}

\vspace{-7pt}

\NI{}can be interpreted as the application of an idealized Runge--Kutta 
method with stages being indexed by an interval \INT01{} (Plato's form 
of a Runge--Kutta method, it is called ``the Picard method'' in 
\citep[p.~97]{But72}, \citep[pp.~153 and 174]{But21}). In the upper half 
of the following table:

\vspace{-14pt}

 \begin{gather*}
   \begin{tabular}{r|l}
     $\VEC{1}$ & function $1(\theta) \equiv 1$ \\
     $\MAT{A} \STRUT$ & operator $u(\theta) \mapsto \int_0^\theta \d 
\theta' \DWS u(\theta')$ \\
     $\VEC{b}$ & functional $u(\theta) \mapsto \int_0^1 \d \theta \DWS
u(\theta)$ \\
     $\VEC{u} . \VEC{v} \vphantom{\vert_{\vert_{\vert_\vert}}} \STRUT$ & 
point-wise product $u(\theta) \mkern0.5mu v(\theta)$ \\
     \hline
     $\VEC{c} \vphantom{\vert^{\vert^{\vert^{-}}}}$ & $\theta$ \\
     $\VEC{q}_n \vphantom{\vert_\vert} $ & $0$ \\
     $\VEC{\Phi}(\textrm{t}) \STRUT$ & $\vert \textrm{t} \vert \mkern2mu 
\theta^{\mkern1mu \vert \textrm{t} \vert - 1} \mkern-2mu / \mkern1mu 
\textrm{t}!$ \\
     $\BPHI(\textrm{t}) \vphantom{\vert^{\vert^{\vert}}}$ & $1 / 
\mkern1mu \textrm{t}!$
   \end{tabular} \end{gather*}

\vspace{-6pt}

\NI{}the process$\mkern1.25mu$/$\mkern0.5mu$result of the substitution 
of entries in the left column with what is on the right will be called a 
\emph{transfiguration}. It can be viewed in two ways: 1) considering 
quantities or statements in the idealized method; and 2) checking for 
desired properties of a Runge--Kutta method, such as simplifying 
assumptions, see, \EG, \citep[p.~52]{But64a}, \citep[s.~321]{But16}, 
\citep[pp.~175, 182, and 208]{HNW93}. Here $\vert \textrm{t} \vert$ is 
the order of tree $\textrm{t}$, \IE, the number of vertices in 
$\textrm{t}$. The factorial $\textrm{t}!$ is recursively defined as 
$\begin{picture}(4.2,6)(0,0) \put(3.35,2.8){\makebox(0,0){\circle*{3}}} 
\end{picture} ! = 1$ and if $\mathrm{t} = [ \mkern1mu \textrm{t}_1 
\mkern2mu \textrm{t}_2 \mkern2mu \DOTS \mkern2mu \textrm{t}_n \mkern1mu 
]$, then $\mathrm{t}! = \vert \textrm{t} \vert \mkern1mu 
\prod_{\mkern1mu m = 1}^{\mkern1mu n} \bigl( \textrm{t}_m \bigr) 
\mkern-0.5mu ! \mkern1mu$.

A Runge--Kutta method of order $4$ should satisfy the order conditions 
$\BPHI(\textrm{t}) = 1 / \textrm{t}!$ for all rooted trees $\textrm{t}$ 
with $\vert \textrm{t} \vert \le 4$:

\vspace{-14pt}

 \begin{gather*}
   \begin{picture}(312,150)(0,0)
     \put(16.9,131){\LEAF}
     \put(18.55,107){\ROOT}
     \put(24,130){\makebox(0,0){\circle*{3}}}
     \put(40,130){\makebox(0,0){\scalebox{1.4}{$=$}}}
     \put(49,131.4){\makebox(0,0)[l]{\scalebox{1.25}{$1$}}}
     \put(94.9,131){\LEAF}
     \put(96.55,90){\ROOT}
     \put(102,121.5){\makebox(0,0){\rotatebox{90}{\rule[3pt]{17pt}{1.44pt}}}}
     \put(102,130){\makebox(0,0){\circle*{3}}}
     \put(102,113){\makebox(0,0){\circle*{3}}}
     \put(118,113){\makebox(0,0){\scalebox{1.4}{$=$}}}
     \put(127,114){\makebox(0,0)[l]{$\DS\frac{1}{2}$}}
     \put(165.2,126.5){\rotatebox{22.5}{\LEAF}}
     \put(190.2,130.8){\rotatebox{-22.5}{\LEAF}}
     \put(182.55,90){\ROOT}
     \put(183.75,120.35){\makebox(0,0){\rotatebox{120}{\rule[2.78pt]{17pt}{1.44pt}}}}
     \put(192.25,120.35){\makebox(0,0){\rotatebox{60}{\rule[2.6pt]{17pt}{1.44pt}}}}
     \put(179.5,127.7){\makebox(0,0){\circle*{3}}}
     \put(196.5,127.7){\makebox(0,0){\circle*{3}}}
     \put(188,113){\makebox(0,0){\circle*{3}}}
     \put(208,113){\makebox(0,0){\scalebox{1.4}{$=$}}}
     \put(217,114){\makebox(0,0)[l]{$\DS\frac{1}{3}$}}
     \put(257.9,131){\LEAF}
     \put(259.55,73){\ROOT}
     \put(265,121.5){\makebox(0,0){\rotatebox{90}{\rule[3pt]{17pt}{1.44pt}}}}
     \put(265,104.5){\makebox(0,0){\rotatebox{90}{\rule[3pt]{17pt}{1.44pt}}}}
     \put(265,130){\makebox(0,0){\circle*{3}}}
     \put(265,113){\makebox(0,0){\circle*{3}}}
     \put(265,96){\makebox(0,0){\circle*{3}}}
     \put(281,96){\makebox(0,0){\scalebox{1.4}{$=$}}}
     \put(290,97){\makebox(0,0)[l]{$\DS\frac{1}{6}$}}
     \put(0.9,75){\LEAF}
     \put(2.55,0){\ROOT}
     \put(8,65.5){\makebox(0,0){\rotatebox{90}{\rule[3pt]{17pt}{1.44pt}}}}
     \put(8,48.5){\makebox(0,0){\rotatebox{90}{\rule[3pt]{17pt}{1.44pt}}}}
     \put(8,31.5){\makebox(0,0){\rotatebox{90}{\rule[3pt]{17pt}{1.44pt}}}}
     \put(8,74){\makebox(0,0){\circle*{3}}}
     \put(8,57){\makebox(0,0){\circle*{3}}}
     \put(8,40){\makebox(0,0){\circle*{3}}}
     \put(8,23){\makebox(0,0){\circle*{3}}}
     \put(24,23){\makebox(0,0){\scalebox{1.4}{$=$}}}
     \put(33,24){\makebox(0,0)[l]{$\DS\frac{1}{24}$}}
     \put(72.2,53.5){\rotatebox{22.5}{\LEAF}}
     \put(97.2,57.8){\rotatebox{-22.5}{\LEAF}}
     \put(89.55,0){\ROOT}
     \put(95,31.5){\makebox(0,0){\rotatebox{90}{\rule[3pt]{17pt}{1.44pt}}}}
     \put(90.75,47.35){\makebox(0,0){\rotatebox{120}{\rule[2.78pt]{17pt}{1.44pt}}}}
     \put(99.25,47.35){\makebox(0,0){\rotatebox{60}{\rule[2.6pt]{17pt}{1.44pt}}}}
     \put(86.5,54.7){\makebox(0,0){\circle*{3}}}
     \put(103.5,54.7){\makebox(0,0){\circle*{3}}}
     \put(95,40){\makebox(0,0){\circle*{3}}}
     \put(95,23){\makebox(0,0){\circle*{3}}}
     \put(114,23){\makebox(0,0){\scalebox{1.4}{$=$}}}
     \put(123,24){\makebox(0,0)[l]{$\DS\frac{1}{12}$}}
     \put(162.2,36.5){\rotatebox{22.5}{\LEAF}}
     \put(186.4,55.7){\LEAF}
     \put(179.55,0){\ROOT}
     \put(193.5,46.2){\makebox(0,0){\rotatebox{90}{\rule[3pt]{17pt}{1.44pt}}}}
     \put(180.75,30.35){\makebox(0,0){\rotatebox{120}{\rule[2.78pt]{17pt}{1.44pt}}}}
     \put(189.25,30.35){\makebox(0,0){\rotatebox{60}{\rule[2.6pt]{17pt}{1.44pt}}}}
     \put(193.5,54.7){\makebox(0,0){\circle*{3}}}
     \put(176.5,37.7){\makebox(0,0){\circle*{3}}}
     \put(193.5,37.7){\makebox(0,0){\circle*{3}}}
     \put(185,23){\makebox(0,0){\circle*{3}}}
     \put(205,23){\makebox(0,0){\scalebox{1.4}{$=$}}}
     \put(214,24){\makebox(0,0)[l]{$\DS\frac{1}{8}$}}
     \put(242.7,35.8){\rotatebox{22.5}{\LEAF}}
     \put(263.9,38){\LEAF}
     \put(278.7,40.1){\rotatebox{-22.5}{\LEAF}}
     \put(265.55,0){\ROOT}
     \put(264,30){\makebox(0,0){\rotatebox{135}{\rule[2.12pt]{20pt}{1.44pt}}}}
     \put(271,30){\makebox(0,0){\rotatebox{90}{\rule[3pt]{14pt}{1.44pt}}}}
     \put(278,30){\makebox(0,0){\rotatebox{45}{\rule[2.12pt]{20pt}{1.44pt}}}}
     \put(257,37){\makebox(0,0){\circle*{3}}}
     \put(271,37){\makebox(0,0){\circle*{3}}}
     \put(285,37){\makebox(0,0){\circle*{3}}}
     \put(271,23){\makebox(0,0){\circle*{3}}}
     \put(296,23){\makebox(0,0){\scalebox{1.4}{$=$}}}
     \put(305,24){\makebox(0,0)[l]{$\DS\frac{1}{4}$}}
   \end{picture} \end{gather*}

\vspace{-7pt}

\NI{}or $\VEC{b} \mkern0.5mu \VEC{1} = 1$, $\VEC{b} \mkern1mu \VEC{c} = 
\frac12$, $\VEC{b} \mkern1mu \VEC{c}^2 = \frac13$, $\VEC{b} \mkern2mu 
\MAT{A} \mkern0.5mu \VEC{c} = \frac{1}{6}$, $\VEC{b} \mkern1mu \VEC{c}^3 
= \frac{1}{4}$, $\VEC{b} \bigl( \VEC{c} . (\MAT{A} \mkern0.5mu \VEC{c}) 
\bigr) = \frac{1}{8}$, $\VEC{b} \mkern2mu \MAT{A} \mkern0.5mu \VEC{c}^2 
= \frac{1}{12}$, $\VEC{b} \mkern2mu \MAT{A}^2 \VEC{c} = \frac{1}{24}$. 
Consider the diagram $\VEC{b} \mkern2mu \MAT{A} \mkern0.5mu \VEC{c}$ in 
the top right statement $\VEC{b} \mkern2mu \MAT{A} \mkern0.5mu \VEC{c} = 
\frac{1}{6}$. If the root or the leaf (which are there out of ``all 
parts of the analytical expression a diagram corresponds to are 
explicitly drawn'' sentiment) are removed, then the diagram becomes the 
vector $\MAT{A} \mkern0.5mu \VEC{c} = \MAT{A}^2 \VEC{1}$ or the row 
vector $\VEC{b} \mkern2mu \MAT{A}^2$, respectively. If both the root and 
the leaf are removed, then the diagram becomes the matrix $\MAT{A}^2$. 
The transfiguration of, \EG, the statement $\VEC{b} \mkern2mu \MAT{A} 
\mkern0.5mu \VEC{c}^2 = \frac{1}{12}$ is

\vspace{-12.5pt}

 \begin{align*}
  & \int\limits_0^1 \d \theta \int\limits_0^\theta \d \theta' \Biggl( 
\int\limits_0^{\theta'} \d \theta'' \Biggr) \Biggl( 
\int\limits_0^{\theta'} \d \theta''' \Biggr) = \int\limits_0^1 \d \theta 
\int\limits_0^\theta \d \theta' \DWS \theta'^2 = \int\limits_0^1 \d 
\theta \DWS \frac{\theta^3}{3} = \frac{1}{12} \\
  & \mkern20mu \VEC{b} \mkern33mu \MAT{A} \mkern15mu \bigl( \mkern9mu 
(\MAT{A} \VEC{1}) \mkern19mu . \mkern19mu (\MAT{A} \VEC{1}) \mkern12mu 
\bigr) \mkern4mu = \mkern17mu \VEC{b} \mkern31mu \MAT{A} \mkern22mu 
\VEC{c}^2
 \end{align*}

\vspace{-4.5pt}

Let $D \bigl( \VEC{\Phi}(\textrm{t}_1) \SEP \VEC{\Phi}(\textrm{t}_2) 
\bigr)$ denote the following property (see \citep[eq.~(2.5)]{AuCh98a}):

\vspace{-17pt}

 \begin{gather*}
   \begin{picture}(260,81)(-2,0)
     \put(-0.65,18){\rotatebox{30}{\BIRCH}}
     \put(34,39.7){\rotatebox{-22.5}{\LARCH}}
     \put(29.55,0){\ROOT}
     \put(40.28,29.72){\makebox(0,0){\rotatebox{45}{\rule{17pt}{1.44pt}}}}
     \put(47,35){\makebox(0,0){\circle*{3}}}
     \put(35,23){\makebox(0,0){\circle*{3}}}
     \put(-0.8,59){\makebox(0,0){$\textrm{t}_1$}}
     \put(46.5,77){\makebox(0,0){$\textrm{t}_2$}}
     \put(85,22.5){\makebox(0,0){\scalebox{2}{$+$}}}
     \put(92.85,31.35){\rotatebox{22.5}{\BIRCH}}
     \put(122.7,29.3){\rotatebox{-30}{\LARCH}}
     \put(129.55,0){\ROOT}
     \put(126.72,29.72){\makebox(0,0){\rotatebox{-45}{\rule{17pt}{1.44pt}}}}
     \put(123,35){\makebox(0,0){\circle*{3}}}
     \put(135,23){\makebox(0,0){\circle*{3}}}
     \put(182,21.5){\makebox(0,0){\scalebox{2}{$=$}}}
     \put(204.55,23.45){\BIRCH}
     \put(211.55,0){\ROOT}
     \put(217,23){\makebox(0,0){\circle*{3}}}
     \put(230,23){\LARCH}
     \put(238.55,0){\ROOT}
     \put(244,23){\makebox(0,0){\circle*{3}}}
   \end{picture} \\
   \VEC{b} \bigl(\VEC{\Phi}(\textrm{t}_1) \mkern1mu . \mkern-1mu \bigl( 
\MAT{A} \VEC{\Phi}(\textrm{t}_2) \bigr) \bigr) + \VEC{b} \bigl( \bigl( 
\MAT{A} \VEC{\Phi}(\textrm{t}_1) \bigr) . \mkern1mu 
\VEC{\Phi}(\textrm{t}_2) \bigr) = \bigl( \VEC{b} \mkern2mu 
\VEC{\Phi}(\textrm{t}_1) \bigr) \bigl( \VEC{b} \mkern2mu 
\VEC{\Phi}(\textrm{t}_2) \bigr)
 \end{gather*}

\vspace{-6pt}

\NI{}whose transfiguration or continuous analog $D(u \SEP \mkern-1mu v)$ 
would be

\vspace{-15pt}

 \begin{gather*}
   \int\limits_0^1 \d \theta \DWS u(\theta) \int\limits_0^\theta \d 
\theta' \DWS v(\theta') \; + \, \int\limits_0^1 \d \theta \Biggl( 
\int\limits_0^\theta \d \theta' \DWS u(\theta') \Biggr) v(\theta) \, = 
\, \Biggl( \int\limits_0^1 \d \theta \DWS u(\theta) \Biggr) \Biggl( 
\int\limits_0^1 \d \theta \DWS v(\theta) \Biggr) \\[1pt]
  \includegraphics{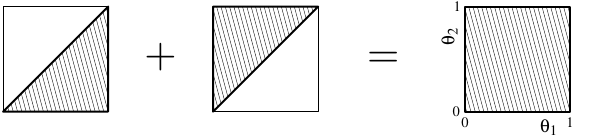} \end{gather*}

\vspace{-12pt}

\NI{}where, assuming it is $u(\theta_1) \mkern0.5mu v(\theta_2)$ being 
integrated, the shading shows the regions of integration in the 
$(\theta_1 \SEP \theta_2)$-plane. Another way to interpret $D(u \SEP 
\mkern-1mu v)$ is as the integration by parts formula. Consider the 
following bilinear form $M$:

\vspace{-13.5pt}

 \begin{gather*}
   \begin{picture}(240,37)(0,0)
     \qbezier(1.5,35)(13.5,3)(25.5,35)
     \put(20.5,18){\makebox(0,0){\color{white}\circle*{14}}}
     \put(13.6,18){\makebox(0,0){$M$}}
     \put(20.5,18){\makebox(0,0){\circle{14}}}
     \put(3,35){\makebox(0,0){\circle*{3}}}
     \put(27,35){\makebox(0,0){\circle*{3}}}
     \put(48,29){\makebox(0,0){\scalebox{1.7}{$=$}}}
     \put(79.55,0){\ROOT}
     \put(76.72,29.72){\makebox(0,0){\rotatebox{135}{\rule{17pt}{1.44pt}}}}
     \put(76.72,29.72){\makebox(0,0){\rotatebox{135}{\color{white}\rule{17pt}{0.72pt}}}}
     \put(90.28,29.72){\makebox(0,0){\rotatebox{45}{\rule{17pt}{1.44pt}}}}
     \put(73,35){\makebox(0,0){\circle*{3}}}
     \put(97,35){\makebox(0,0){\circle*{3}}}
     \put(85,23){\makebox(0,0){\circle*{3}}}
     \put(119,29){\makebox(0,0){\scalebox{1.7}{$+$}}}
     \put(147.55,0){\ROOT}
     \put(144.72,29.72){\makebox(0,0){\rotatebox{135}{\rule{17pt}{1.44pt}}}}
     \put(158.28,29.72){\makebox(0,0){\rotatebox{45}{\rule{17pt}{1.44pt}}}}
     \put(158.28,29.72){\makebox(0,0){\rotatebox{45}{\color{white}\rule{17pt}{0.72pt}}}}
     \put(141,35){\makebox(0,0){\circle*{3}}}
     \put(165,35){\makebox(0,0){\circle*{3}}}
     \put(153,23){\makebox(0,0){\circle*{3}}}
     \put(186,29){\makebox(0,0){\scalebox{1.7}{$-$}}}
     \put(207.55,12){\ROOT}
     \put(213,35){\makebox(0,0){\circle*{3}}}
     \put(231.55,12){\ROOT}
     \put(237,35){\makebox(0,0){\circle*{3}}}
   \end{picture} \end{gather*}

\vspace{-8.5pt}

\NI{}whose matrix $\MAT{M} = \bigl[ m_{\mkern1mu i \mkern-1mu j} 
\bigr]$, with $m_{\mkern1mu i \mkern-1mu j} = b_{i \mkern1mu} \AIJ + 
b_{\mkern-1.2mu j \mkern1.2mu} a_{\mkern-1.2mu j \mkern1.2mu i} - b_{i 
\mkern1mu} b_{\mkern-1mu j \mkern1mu}$, is symmetric. The property 
$D(\VEC{u} \SEP \mkern-1mu \VEC{v})$ can be written as $M(\VEC{u} \SEP 
\mkern-1mu \VEC{v}) = \VEC{u}^{\textrm{T}} \MAT{M} \mkern1mu \VEC{v} = 
0$.

An order $p$ Runge--Kutta method is said to be pseudo-symplectic of 
order $(p \SEP q)$ if $D \bigl( \VEC{\Phi}(\textrm{t}_1) \SEP 
\VEC{\Phi}(\textrm{t}_2) \bigr)$ for all trees $\mathrm{t}_1$ and 
$\mathrm{t}_2$ such that $\vert \mathrm{t}_1 \vert + \vert \mathrm{t}_2 
\vert \le q$ \citep[cor.~2.2 and eq.~(2.7)]{AuCh98a}. $B$-series 
coefficients, $\mathrm{t} \mapsto \BPHI(\textrm{t})$, for a method of 
order $(p \SEP q)$ are in the coset $\bigl( \mathrm{E} + \mathsf{O}_{p + 
1} \bigr)$ and in the subgroup $\bigl( \cap_{\mkern1mu r = 1}^{\mkern1mu 
q - 1} \mathsf{D}_{\mkern1.25mu r \SEP q - r} \bigr)$ of the Butcher 
group $\mathsf{B}$ \citep[pp.~166 and 173]{But21}.

The property $\bigl( D \bigl( \VEC{\Phi} \SEP \VEC{u} \bigr)$ for all 
$\VEC{u} \in \pmb{\mathbf{R}}^s \bigr)$ will be shortened as $\VEC{D} 
\bigl( \VEC{\Phi} \bigr)$ (see \citep[def.~3.1]{AuCh98a}). It is 
equivalent to the linear form $\VEC{u} \mapsto M \bigl( \VEC{\Phi} \SEP 
\VEC{u} \bigr)$ being a zero function, or $\MAT{M} \mkern0.5mu 
\VEC{\Phi} = \VEC{0}$. The \NTH{\smash{j{\mkern1.2mu}}} component of 
$\VEC{D} \bigl( \VEC{\Phi} \bigr)$, the property $D \bigl( \VEC{\Phi} 
\SEP \VEC{e}_{\mkern-1.5mu j} \bigr)$, where $\VEC{e}_{\mkern-1.5mu j}$ 
is the \NTH{\smash{j{\mkern1.2mu}}} vector of the standard basis of 
$\pmb{\mathbf{R}}^s$, will be denoted by $D_{\mkern-1.5mu j} \bigl( 
\VEC{\Phi} \bigr)$. Assuming that $\VEC{b} \mkern0.5mu \VEC{1} = 1$ and 
$\VEC{b} \mkern1mu \VEC{c} = \smash{\frac{1}{2}}$, here are the diagrams 
for $\VEC{D}\BRA1$ and $\VEC{D}(\VEC{c})$:

\vspace{-15pt}

 \begin{gather*}
   \begin{picture}(318,162)(0,0)
     \put(130.55,113.45){\BIRCH}
     \put(137.53,75){\ROOT}
     \put(143,104.5){\makebox(0,0){\rotatebox{90}{\rule[2.88pt]{17pt}{1.44pt}}}}
     \put(143,113){\makebox(0,0){\circle*{3}}}
     \put(143,96){\makebox(0,0){\circle*{3}}}
     \put(172,95.5){\makebox(0,0){\scalebox{2}{$+$}}}
     \put(187.7,106.8){\rotatebox{22.5}{\LEAF}}
     \put(203.3,101.9){\rotatebox{-30}{\BIRCH}}
     \put(208.55,73){\ROOT}
     \put(205.72,102.72){\makebox(0,0){\rotatebox{-45}{\rule{17pt}{1.44pt}}}}
     \put(202,108){\makebox(0,0){\circle*{3}}}
     \put(214,96){\makebox(0,0){\circle*{3}}}
     \put(255,94.5){\makebox(0,0){\scalebox{2}{$=$}}}
     \put(273.55,96.45){\BIRCH}
     \put(280.55,73){\ROOT}
     \put(286,96){\makebox(0,0){\circle*{3}}}
     \put(32,113){\makebox(0,0)[lb]{$\DS {\VEC{b} \mkern1.5mu \MAT{A} 
\mkern0.5mu \VEC{u}} + {\VEC{b} \mkern0.5mu (\VEC{c} . \mkern1mu 
\VEC{u})} = {\VEC{b} \mkern1mu \VEC{u}}$}}
     \put(-3.3,33.8){\rotatebox{22.5}{\LEAF}}
     \put(17.5,37.7){\LARCH}
     \put(17.55,0){\ROOT} 
     \put(14.72,29.72){\makebox(0,0){\rotatebox{-45}{\rule{17pt}{1.44pt}}}} 
     \put(27.25,30.35){\makebox(0,0){\rotatebox{60}{\rule[2.6pt]{17pt}{1.44pt}}}} 
     \put(11,35){\makebox(0,0){\circle*{3}}}
     \put(31.5,37.7){\makebox(0,0){\circle*{3}}} 
     \put(23,23){\makebox(0,0){\circle*{3}}} 
     \put(55,22.5){\makebox(0,0){\scalebox{1.8}{$+$}}} 
     \put(69.5,52.4){\LEAF}
     \put(80.7,29.3){\rotatebox{-30}{\LARCH}}
     \put(87.55,0){\ROOT} 
     \put(78.54,44){\makebox(0,0){\rotatebox{105}{\rule[2.78pt]{17pt}{1.44pt}}}}
     \put(84.72,29.72){\makebox(0,0){\rotatebox{-45}{\rule{17pt}{1.44pt}}}} 
     \put(76.6,51.4){\makebox(0,0){\circle*{3}}}
     \put(81,35){\makebox(0,0){\circle*{3}}} 
     \put(93,23){\makebox(0,0){\circle*{3}}} 
     \put(131,21.5){\makebox(0,0){\scalebox{1.8}{$=$}}} 
     \put(147,23){\makebox(0,0)[l]{$\DS\frac{1}{2}$}}
     \put(155,23){\LARCH}
     \put(163.55,0){\ROOT} 
     \put(169,23){\makebox(0,0){\circle*{3}}} 
     \put(318,40){\makebox(0,0)[rb]{$\DS {\VEC{b} \bigl( \VEC{c} . (\MAT{A} 
\mkern0.5mu \VEC{u}) \mkern-1mu \bigr)} + {\VEC{b} \bigl( \mkern-1mu 
(\MAT{A} \mkern0.5mu \VEC{c}) . \mkern1mu \VEC{u} \bigr)} = 
{{\TS\frac12} \mkern0.5mu \VEC{b} \mkern1mu \VEC{u}}$}} \end{picture} 
\end{gather*}

\vspace{-9pt}

\NI{}The property $\VEC{D}\BRA1$ coincides with $D(1)$ 
\citep[p.~52]{But64a}, \citep[fig.~321(ii)]{But16}, \citep[pp.~175 and 
208]{HNW93}, \citep[p.~193]{But21}. If $C(2)$ and $D(1)$ are satisfied, 
then $\VEC{D}(\VEC{c})$ is equivalent to $D(2)$. For an explicit 
Runge--Kutta method that satisfies the explicit version $\bigl( (\MAT{A} 
\mkern1mu \VEC{c})_i \ne c_i^2 / 2 \bigr) \Rightarrow \bigl( (i = 2) 
\land (b_2 = 0) \bigr)$ of the property $C(2)$, it is impossible to have 
both $\VEC{D}\BRA1$ and $\VEC{D}(\VEC{c})$ \citep[prop.~3.6]{AuCh98a}.

Let a function $u(\theta)$ be called \INT01{}-even or \INT01{}-odd, if 
$u(1 - \theta) = u(\theta)$ or $u(1 - \theta) = -u(\theta)$, 
respectively. The following parity properties are true:

\vspace{-14.5pt}

 \begin{gather} \begin{array}{l}
   \!\!\!\!{\DS \mbox{if $u(\theta)$ is \INT01{}-odd, then $U(\theta) = 
\smash{\int_0^\theta \d \theta' \, u(\theta')}$ is \INT01{}-even}} 
\\[6pt]
   \!\!\!\!{\DS \mbox{if $u(\theta)$ is \INT01{}-even and $\smash{\int_0^1 \d 
\theta \, u(\theta) = 0}$, then $U(\theta) = \smash{\int_0^\theta \d 
\theta' \, u(\theta')}$ is \INT01{}-odd}}\!\!\!
  \end{array} \label{even_odd} \end{gather}

\vspace{-8.5pt}

\section{Family of pseudo-symplectic methods} \label{family}

\hskip\parindent{}A pseudo-symplectic Runge--Kutta method of order $(4 
\SEP 8)$ will be searched for within $8$-stage explicit methods with 
$c_4 = c_5 = \frac12$, $c_6 = 1 - c_3$, $c_7 = 1 - c_2$, $c_8 = 1$, $b_6 
= b_3$, $b_7 = b_2$, and $b_8 = b_1$. Vectors $\VEC{x}$ and $\VEC{y}$ of 
$8$ components are going to be called ``even'' and ``odd'' if $x_6 = 
x_3$, $x_7 = x_2$, $x_8 = x_1$, and $y_1 + y_8 = y_2 + y_7 = y_3 + y_6 = 
y_4 = y_5 = 0$, respectively. Let $\VEC{p}_1 = 2 \VEC{c} - \VEC{1}$, 
$\VEC{p}_2 = 6 \mkern1.5mu \VEC{c}^2 - 6 \mkern1.5mu \VEC{c} + \VEC{1}$, 
and $\VEC{p}_3 = 20 \mkern1.5mu \VEC{c}^3 - 30 \mkern1.5mu \VEC{c}^2 + 
12 \VEC{c} - \VEC{1}$ be the analogs of shifted Legendre polynonials, 
orthogonal on \INT01{} (see, \EG, \citep[s.~342]{But16}, 
\citep[p.~203]{But21}). The vectors $\VEC{1}$, $\VEC{b}^{\textrm{T}}$, 
and $\VEC{p}_2$ are ``even''; while $\VEC{p}_1$ and $\VEC{p}_3$ are 
``odd''. With such a definition of ``even'' and ``odd'' vectors, the 
eq.~(\ref{even_odd}) could be viewed as the transfiguration of the 
statements

\vspace{-13pt}

 \begin{gather*} \begin{array}{l}
   {\DS \mbox{if $\VEC{u}$ is ``odd'', then $\MAT{A} \mkern0.5mu 
\VEC{u}$ is ``even''}} \\[2pt]
   {\DS \mbox{if $\VEC{u}$ is ``even'' and $\VEC{b} \mkern1mu \VEC{u} = 
0$, then $\MAT{A} \mkern0.5mu \VEC{u}$ is ``odd''}}
  \end{array} \end{gather*}

\vspace{-4pt}

\NI{}To imitate these parity properties, it will be demanded that 
$\MAT{A} \VEC{p}_1$ is ``even'' (this implies that the vector $\VEC{q}_1 
= \MAT{A} \VEC{c} - \frac12 \VEC{c}^2 = \frac12 \MAT{A} \VEC{p}_1 + 
\frac12 \VEC{c} (\VEC{1} - \VEC{c})$ is ``even'' too), and that $\MAT{A} 
\VEC{p}_2$ and $\MAT{A} \mkern0.5mu \VEC{q}_1$ are ``odd''. These 
demands can be viewed as simplifying assumptions, as they increase the 
redundancy in the order conditions (for example, the statements $\VEC{b} 
\VEC{1} = 1$ and $\VEC{b} \VEC{p}_2 = \VEC{b} \mkern1mu \VEC{q}_1 = 0$ 
imply that a method is of order $4$). Also the properties $\VEC{D}\BRA1$ 
and $\VEC{D}(\VEC{c})$ (or an easier to exploit $\VEC{D}\BRAP1$) will be 
assumed.\footnote{~The derivation was done in interaction with computer 
algebra system Wolfram Mathematica~8.0, mainly using commands 
\scalebox{0.8}[0.825]{\pmb{\scriptsize\tt\bfseries Solve}} to 
symbolically solve linear equations, 
\scalebox{0.8}[0.825]{\pmb{\scriptsize\tt Simplify}}, and 
\scalebox{0.8}[0.825]{\pmb{\scriptsize\tt\bfseries Factor}}. 
Step-by-step description of the process (of course, this could be done 
in many ways) allows for an easy reproduction (also see Appendix).}

The vectors $\VEC{p}_3 - \VEC{p}_1$, $\VEC{p}_1 \mkern-0.5mu . 
\mkern1.75mu \VEC{q}_1$, $\MAT{A} \VEC{p}_2$, and $\MAT{A} \mkern0.5mu 
\VEC{q}_1$ are ``odd'' and have the \FIRST{} and the \NTH{8} components 
being equal to $0$. Thus, they lie in the span of two vectors: 
$\VEC{v}_2 = \bigl[ \; 0 ~~~ 1 ~~~ 0 ~~~ 0 ~~~ 0 ~~~ 0 ~~ -\!1 ~~~ 0 \; 
\bigr]{}^{\textrm{T}}$ and $\VEC{v}_3 = \bigl[ \; 0 ~~~ 0 ~~~ 1 ~~~ 0 
~~~ 0 ~~ -\!1 ~~~ 0 ~~~ 0 \; \bigr]{}^{\textrm{T}}$. To satisfy the 
order conditions $D \bigl( \VEC{\Phi}(\textrm{t}_1), 
\VEC{\Phi}(\textrm{t}_2) \bigr)$ with $\vert \textrm{t}_1 \vert = \vert 
\textrm{t}_2 \vert = 4$, it is necessary and sufficient that $\VEC{b} 
\bigl( \VEC{v}_2 . (\MAT{A} \VEC{v}_2) \bigr) = -a_{72} b_2 = 0$, 
$\VEC{b} \bigl( \VEC{v}_3 . (\MAT{A} \VEC{v}_3) \bigr) = -a_{63} b_3 = 
0$, and $\VEC{b} \bigl( \VEC{v}_2 . (\MAT{A} \VEC{v}_3) \bigr) + \VEC{b} 
\bigl( (\MAT{A} \VEC{v}_2) . \VEC{v}_3 \bigr) = (a_{76} - a_{73}) b_2 + 
(a_{32} - a_{62}) b_3 = 0$. To support the observance of the parity 
properties and not to lose the flexibility in choosing $b_2$ and $b_3$, 
let $a_{63} = a_{72} = 0$, $a_{62} = a_{32}$, and $a_{76} = a_{73}$.

The conditions $-c_2 = \BRAAP1_2 = \BRAAP1_7 = -a_{71}$ and $\BRAAP1_3 = 
\BRAAP1_6$ imply $a_{71} = c_2$ and $a_{61} = a_{31}$, respectively. The 
property $D_1{\BRAP1}$ becomes $a_{81} b_1 = 0$, thus $a_{81} = 0$.

The property $D_2{\BRAP1}$ implies $a_{82} = c_2 b_2 / b_1$. Both 
$D_7{\BRA1}$ and $D_7{\BRAP1}$ imply $a_{87} = c_2 b_2 / b_1$. Now all 
three $D_3{\BRAP1} \Leftrightarrow D_6{\BRAP1}$, $D_8{\BRAP1}$, and 
$\BRAAP1_8 = 0$ imply $a_{86} = a_{83}$.

If $b_2 = 0$, then $D_6{\BRA1}$ implies $a_{83} = c_3 b_3 / b_1$, and 
then $D_6{\BRAP1}$ means $b_3 a_{32} c_2 = 0$. Both cases $b_3 = 0$ and 
$a_{32} = 0$ lead to contradictions, thus $b_2 \ne 0$. With $a_{31} = 
c_3 - a_{32}$, the properties $D_6{\BRA1}$ and $D_6{\BRAP1}$ imply 
$a_{83} = a_{31} b_3 / b_1$ and $a_{73} = a_{32} b_3 / b_2$. The Butcher 
tableau \citep[p.~191]{But64b} now looks like

\vspace{-14.1pt}

 \begin{gather*}
   \begin{array}{r|ccccccccc}
     0\;\,\vphantom{\vert_\vert} \\
     {c_2}\;\,\STRUT & \phantom{-}c_2 \\
     {c_3}\;\,\STRUT & \phantom{-}a_{31} & \phantom{-}a_{32} \\
     \frac{1}{2}\;\,\STRUT & \phantom{-}a_{41} & \phantom{-}a_{42} & 
\phantom{-}a_{43} \\
     \frac{1}{2}\;\,\STRUT & \phantom{-}a_{51} & \phantom{-}a_{52} & 
\phantom{-}a_{53} & \phantom{-}a_{54} \\
     {1 - c_3}\;\,\STRUT & \,\phantom{-}a_{31}\, & \phantom{-}a_{32} & 
\phantom{-}0 & \,\phantom{-}a_{64}\, & \,\phantom{-}a_{65}\, \\
     {1 - c_2}\;\,\STRUT & \phantom{-}c_2 & \phantom{-}0 & 
\phantom{-}a_{32} \frac{b_3}{b_2} & \phantom{-}a_{74} & 
\phantom{-}a_{75} & \phantom{-}a_{32} \frac{b_3}{b_2} \\
     1\;\,\STRUT\vphantom{\vert_{\vert_\vert}} & \phantom{-}0 & 
\phantom{-}c_{2} \frac{b_2}{b_1} & \phantom{-}a_{31} \frac{b_3}{b_1} & 
\phantom{-}a_{84} & \phantom{-}a_{85} & \phantom{-}a_{31} 
\frac{b_3}{b_1} & \phantom{-}c_{2} \frac{b_2}{b_1} \\
     \hline
     \vphantom{\vert^{\vert^\vert}} & \phantom{-}b_1 & \phantom{-}b_2 & 
\phantom{-}b_3 & \phantom{-}b_4 & \phantom{-}b_5 & \phantom{-}b_3 & 
\phantom{-}b_2 & \,\phantom{-}b_1
   \end{array} \end{gather*} 

\vspace{-5.1pt}

As ${\BRAA1}_6 = 1 - c_3$, let $a_{64} = 1 - 2 c_3 - a_{65}$. Then 
$(\MAT{A} \VEC{p}_2)_3 + (\MAT{A} \VEC{p}_2)_6 = 0$ implies $c_3 \ne 
\frac{1}{6}$ and $a_{32} = (6 c_3 - 1) / 24 c_2 (1 - c_2)$. As 
${\BRAA1}_4 = \frac12$, let $a_{41} = \frac12 - a_{42} - a_{43}$. Then 
$(\MAT{A} \VEC{p}_2)_4 = 0$ implies $a_{42} = \bigl( 1 - 12 a_{43} c_3 
(1 - c_3) \bigr) / 12 c_2 (1 - c_2)$. Now $(\MAT{A} \mkern1mu \VEC{q}_1 
\mkern-1.5mu )_4 = 0$ is equivalent to $a_{43} = c_2 / \bigl( 6 c_3 - 12 
c_3 (c_3 - c_2) - 1 \bigr)$, with the denominator necessarily being 
non-zero. As ${\BRAA1}_5 = \frac12$, let $a_{51} = \frac12 - a_{52} - 
a_{53} - a_{54}$. Then $(\MAT{A} \VEC{p}_2)_5 = 0$ implies $a_{54} = 
\frac13 - 4 a_{52} c_2 (1 - c_2) - 4 a_{53} c_3 (1 - c_3)$.

As $(\MAT{A} \VEC{1})_7 = 1 - c_2$, let $a_{74} = 1 - 2 c_2 - a_{75} - 
b_3 (6 c_3 - 1) / 12 b_2 c_2 (1 - c_2)$. Then $(\MAT{A} \VEC{p}_2)_2 + 
(\MAT{A} \VEC{p}_2)_7 = \bigl( b_3 (6 c_3 - 1) (1 - 2 c_3)^2 - 4 b_2 c_2 
(1 - c_2) (1 - 6 c_2) \bigr) / 8 b_2 c_2 (1 - c_2) = 0$, thus $b_3 = 4 
b_2 c_2 (1 - c_2) (1 - 6 c_2) / (1 - 2 c_3)^2 (6 c_3 - 1)$. The 
coefficient $a_{84}$ is determined from ${\BRAA1}_8 = 1$. Then the 
derivative of $b_1 (\MAT{A} \VEC{p}_2)_8$ with respest to $b_1$ is equal 
to $-\frac{1}{2} \ne 0$, so $(\MAT{A} \VEC{p}_2)_8 = 0$ determines $b_1 
= b_2 \bigl( 1 - 12 c_2 + 24 c_2^2 (1 - c_2) - 6 c_3 + 96 c_2 c_3 - 312 
c_2^2 c_3 + 288 c_2^3 c_3 \bigr) / (6 c_3 - 1)$. Now $\MAT{A} \VEC{p}_2 
= \bigl[ \; 0 ~~~ c_2 ~~~ \frac{1 - 2 c_3}{4} ~~~ 0 ~~~ 0 ~~~ 
{-{\mkern-1mu}\frac{1 - 2 c_3}{4}} ~~ {-{\mkern-1mu}c_2} ~~~ 0 \; 
\bigr]{}^{\textrm{T}}$ is ``odd'', also $\VEC{q}_1$ is ``even''. The 
coefficient $a_{85}$ and the weight $b_4$ are determined from 
$D_5(\VEC{1})$ and $\VEC{b} \mkern0.5mu \VEC{1} = 1$, respectively.

Now $\VEC{b} \mkern1mu \VEC{c}^2 = \frac{1}{4} + 2 b_2 c_2 \bigl( 18 c_3 
(1 - 2 c_2)^2 - 1 - 2 c_2 \bigr) / (6 c_3 - 1)$. Thus, $18 c_3 (1 - 2 
c_2)^2 - 1 - 2 c_2 \ne 0$ and $b_2 = (6 c_3 - 1) / 24 c_2 \bigl( 18 c_3 
(1 - 2 c_2)^2 - 1 - 2 c_2 \bigr)$.

The coefficient $a_{75}$ is determined from $D_5{\BRAP1}$.

The properties $D_2{\BRA1}$, $D_3{\BRA1}$, and $\BRAAQ1_5 = 0$ can be 
written as $b_5 a_{52} + \phi_{2b} b_5 = \phi_{20}$, $b_5 a_{53} + 
\phi_{3b} b_5 = \phi_{30}$, and $\phi_2 a_{52} + \phi_3 a_{53} = 
\phi_0$, respectively, where all the coefficients $\phi_{2b}$, 
$\phi_{20}$, $\phi_{3b}$, $\phi_{30}$, $\phi_2$, $\phi_3$, and $\phi_0$ 
are rational functions of $c_2$ and $c_3$. It can be checked that 
$\phi_2 \phi_{2b} + \phi_3 \phi_{3b} + \phi_0 = 0$. The following 
quantity should be equal to zero: $\phi_2 (b_5 a_{52} + \phi_{2b} b_5 - 
\phi_{20}) + \phi_3 (b_5 a_{53} + \phi_{3b} b_5 - \phi_{30}) = b_5 
(\phi_2 a_{52} + \phi_3 a_{53} + \phi_2 \phi_{2b} + \phi_3 \phi_{3b}) - 
\phi_2 \phi_{20} - \phi_3 \phi_{30} = b_5 \bigl( \phi_0 + (-\phi_0) 
\bigr) - \phi_2 \phi_{20} - \phi_3 \phi_{30} = -(\phi_2 \phi_{20} + 
\phi_3 \phi_{30})$. This can be rewritten as

\vspace{-18pt}

 \begin{gather*}
   \zeta(c_2 \SEP c_3) = \sum_{m, n} \zeta_{\mkern1mu m \mkern0.5mu n} 
\mkern1.5mu c_2^m \mkern0.5mu (2 c_3)^n = 0 \qquad
   \begin{array}{r|rrrrr}
     \mbox{\footnotesize$5$}\; & \phantom{-}9 & -9 \\
     \mbox{\footnotesize$4$}\; & {\;}-21 & {\;}-15 & {\;}-108 & 
\phantom{-}216 \\
     \mbox{\footnotesize$3$}\; & \phantom{-}15 & \phantom{-}75 & 
\phantom{-}294 & {\;}-576 \\
     \mbox{\footnotesize$2$}\; & -1 & -93 & -198 & \phantom{-}396 & 
{\;}\phantom{-}72 \\
     \mbox{\footnotesize$n = 1$}\; & -3 & \phantom{-}53 & \phantom{-}18 
& -132 \\
     \mbox{\footnotesize$n = 0$}\; & \phantom{-}1 & -13 & \phantom{-}20 
\\
     \hline
     \zeta_{\mkern1mu m \mkern0.5mu n}\vphantom{\vert^\vert}\; & 
\mbox{\footnotesize$m = 0$} & \mbox{\footnotesize$m = 1$} & 
\mbox{\footnotesize$\phantom{-}2$}~ & 
\mbox{\footnotesize$\phantom{-}3$}~~ & 
\mbox{\footnotesize$\phantom{-}4$}~
   \end{array} \end{gather*}

\vspace{-3pt}

 \begin{figure}
   \vspace{10pt}
   \centerline{\includegraphics{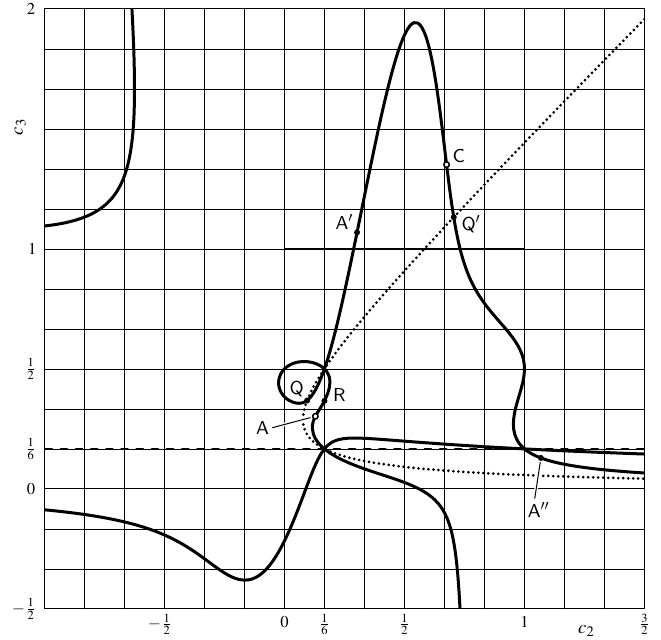}}
   \caption{The curve $\zeta(c_2 \SEP c_3) = 0$ (solid curve); the line 
$c_3 = \frac{1}{6}$ (dashed curve); and the hyperbola $6 c_3 - 12 c_3 
(c_3 - c_2) - 1 = 0$ (dotted curve, the coefficient $a_{43}$ is infinite 
on it) which intersects with the curve $\zeta = 0$ at the points 
$(\frac{1}{6} \SEP \frac{1}{6})$, $(\frac{1}{6} \SEP \frac{1}{2})$, 
$\mathsf{Q} = \bigl( (6 - \sqrt{21}) / 15 \SEP (9 - \sqrt{21}) / 12 
\bigr) = (0.09449..., 0.36811...)$, and $\mathsf{Q}' = \bigl( (6 + 
\sqrt{21}) / 15 \SEP (9 + \sqrt{21}) / 12 \bigr) = (0.70550..., 
1.13188...)$. Some branches of the curve $\zeta = 0$ are not shown as 
they correspond to too large values of $c_2$ or $c_3$ to be visible on 
the graph. The points $\mathsf{C}$ and $\mathsf{A}$ are indicated by 
open dots and correspond to the method in eq.~(\ref{method_C}) and the 
$1$-dimensional family of methods that includes the method in 
eq.~(\ref{the_method}), respectively.}
 \end{figure}

\NI Both coefficients $\phi_{20}$ and $\phi_{30}$ are equal to zero at 
the point $(c_2 \SEP c_3) = \mathsf{C} = (2 \gamma \SEP 4 \gamma)$, 
where $\gamma = \bigl( 2 + 2^{1 / 3} + 2^{-1 / 3} \bigr) / 12 = 1 / 4 
\mkern0.5mu (2 - 2^{1/3}) = 0.33780...$ There $\phi_{2b} = 0$ and 
$\phi_{3b} = 2 \gamma - \frac{1}{2}$, thus $a_{53} = \frac{1}{2} (-5 + 
48 \gamma - 96 \gamma^2) \mkern1mu a_{52} + \frac{1}{2} - 2 \gamma$, 
also $b_5 a_{52} = 0$ (from $D_2{\BRA1}$) and $a_{65} a_{52} = 0$ (from 
$\BRAAQ1_3 + \BRAAQ1_6 = 0$). In both cases of $a_{52} = 0$ (then 
$\VEC{a}_{4*} = \VEC{a}_{5*}$, \IE, the \NTH{4} and the \NTH{5} stages 
duplicate each other) and $a_{52} \ne 0$ (then $b_5 = 0$ and 
$\VEC{a}_{*5} = \VEC{0}$, \IE, the \NTH{5} stage is not used), the 
method becomes equivalent to

\vspace{-12pt}

 \begin{gather}
   \begin{array}{r|cccccccc}
     0\;\,\vphantom{\vert_\vert} \\
     {2 \gamma}\;\,\STRUT & \phantom{-}2 \gamma \\
     {4 \gamma}\;\,\STRUT & \phantom{-}0 & \;\,\phantom{-}4 \gamma \\
     \frac{1}{2}\;\,\STRUT & \phantom{-}2 \gamma & \;\;\phantom{-}0 & 
\phantom{-}\frac{1}{2} - 2 \gamma \\
     {1 - 4 \gamma}\;\,\STRUT & \phantom{-}0 & \;\;\phantom{-}4 \gamma & 
\phantom{-}0 & \phantom{-}1 - 8 \gamma \\
     {1 - 2 \gamma}\;\,\STRUT & \phantom{-}2 \gamma & \;\;\phantom{-}0 & 
\phantom{-}\frac{1}{2} - 2 \gamma & \phantom{-}0 & 
\phantom{-}\frac{1}{2} - 2 \gamma \\
     {1}\;\,\STRUT & \phantom{-}0 & \;\;\phantom{-}4 \gamma & 
\phantom{-}0 & \phantom{-}1 - 8 \gamma & \phantom{-}0 & \,\phantom{-}4 
\gamma \\
     \hline
     \vphantom{\vert^{\vert^{\vert^\vert}}} & \phantom{-}\gamma & 
\;\;\phantom{-}2 \gamma & \phantom{-}\frac{1}{4} - \gamma & 
\phantom{-}\frac{1}{2} - 4 \gamma & \phantom{-}\frac{1}{4} - \gamma & 
\;\phantom{-}2 \gamma & \;\;\,\phantom{-}\gamma
   \end{array} \label{method_C} \end{gather}

\vspace{-3pt}

\NI{}This $7$-stage pseudo-symplectic Runge--Kutta method of order $(4 
\SEP 9)$ is strongly connected to a well-known explicit symplectic 
method for separated Hamiltonians \citep[eqs.~(4.9) and (4.12)]{FoRu90}, 
\citep[eq.~(2.11)]{Yos90}.

Besides the point $\mathsf{C}$ one must have $b_5 \ne 0$, and thus 
$a_{52} = \phi_{20} / b_5 - \phi_{2b}$ and $a_{53} = \phi_{30} / b_5 - 
\phi_{3b}$. The cubic equation $z \mkern1mu (z - \frac{1}{2}) (z - 1) = 
\frac{1}{24}$ has three roots: $z_1 = 1 / 2 - \sin(2 \pi / 9) / \sqrt{3} 
= 0.12888...$, $z_2 = 1 / 2 - \sin(2 \pi / 9 + 2 \pi / 3) / \sqrt{3} = 1 
/ 8 \sin^2 (2 \pi / 9) = 0.30253...$, and $z_3 = 1 / 2 - \sin(2 \pi / 9 
+ 4 \pi / 3) / \sqrt{3} = 1.06857...$ The quantity $b_5 \bigl( \BRAAQ1_3 
+ \BRAAQ1_6 \bigr)$, which should be equal to zero, is a linear function 
of $b_5$ and $a_{65}$, and its derivative with respect to $b_5$ is not 
equal to zero besides the points $(c_2 \SEP c_3) = \mathsf{A} = (z_1 
\SEP z_2)$, $\mathsf{A}' = (z_2 \SEP z_3)$, and $\mathsf{A}'' = (z_3 
\SEP z_1)$. There the weight $b_5$ is determined from $\BRAAQ1_3 + 
\BRAAQ1_6 = 0$, and then $a_{65}$ is found from $\BRAAQ1_8 = 0$. The 
vector $\MAT{A} \mkern0.5mu \VEC{q}_1 = \bigl[ \; 0 ~~~ 0 ~~~ 
{-{\mkern-1mu}\frac{c_2 (6 c_3 - 1)}{48 (1 - c_2)}} ~~~ 0 ~~~ 0 ~~~ 
{\frac{c_2 (6 c_3 - 1)}{48 (1 - c_2)}} ~~ 0 ~~~ 0 \; 
\bigr]{}^{\textrm{T}}$ is now ``odd''. The result is a $1$-dimensional 
(along the curve $\zeta(c_2 \SEP c_3) = 0$) family of pseudo-symplectic 
methods of order $(4 \SEP 6)$ that satisfy $\VEC{D}{\BRA1}$, 
$\VEC{D}(\VEC{c})$, and $\VEC{D}(\VEC{c}^2)$, but not $\VEC{D}(\MAT{A} 
\mkern0.5mu \VEC{c})$. For example, at the point $(c_2 \SEP c_3) = 
\mathsf{R} = (\frac{1}{6} \SEP \frac{11}{30})$ a method is equivalent 
(the \NTH{6} stage is not used) to

\vspace{-13.5pt}

 \begin{gather*}
   \begin{array}{r|cccccccc}
     0\;\,\vphantom{\vert_\vert} \\
     {\frac{1}{6}}\;\,\STRUT & \;\phantom{-}\frac{1}{6} \\
     {\frac{11}{30}}\;\,\STRUT & \;\phantom{-}\frac{1}{150} & 
\;\phantom{-}\frac{9}{25} \\
     \frac{1}{2}\;\,\STRUT & \;\phantom{-}\frac{1}{4} & \;-\frac{13}{48} 
& \;\phantom{-}\frac{25}{48} \\
     \frac{1}{2}\;\,\STRUT & \;-\frac{37}{50} & 
\;\phantom{-}\frac{59}{25} & \;-2 & \;\phantom{-}\frac{22}{25} \\
     {\frac{5}{6}}\;\,\STRUT & \;\phantom{-}\frac{1}{6} & \;\phantom{-}0 
& \;\phantom{-}0 & \;\phantom{-}\frac{4}{11} & 
\;\phantom{-}\frac{10}{33} \\
     1\;\,\STRUT\vphantom{\vert_{\vert_\vert}} & \;\phantom{-}0 & 
\;\phantom{-}\frac{3}{8} & \;\phantom{-}0 & \;\phantom{-}\frac{4}{11} & 
\;-\frac{5}{44} & \;\phantom{-}\frac{3}{8} \\
     \hline
     \vphantom{\vert^{\vert^{\vert^\vert}}} & \;\phantom{-}\frac{1}{12} 
& \;\phantom{-}\frac{3}{16} & \;\phantom{-}0 & \;\phantom{-}\frac{4}{11} 
& \;\phantom{-}\frac{25}{264} & \;\phantom{-}\frac{3}{16} & 
\;\phantom{-}\frac{1}{12}
   \end{array} \end{gather*}

\vspace{-3pt}

At all the three points $\mathsf{A}$, $\mathsf{A}'$, and $\mathsf{A}''$ 
the vector $\MAT{A} \mkern0.5mu \VEC{q}_1$ is ``odd'', plus 
$\VEC{D}{\BRA1}$ and $\VEC{D}(\VEC{c})$ are satisfied. The coefficient 
$a_{65}$ is determined from $D_5(\VEC{c}^2)$. Now $\VEC{D}(\VEC{c}^2)$ 
and $\VEC{D}(\MAT{A} \mkern0.5mu \VEC{c})$ are satisfied, and the result 
is a $1$-dimensional family of $8$-stage pseudo-symplectic Runge--Kutta 
methods of order $(4 \SEP 8)$ indexed by a parameter $\psi$:

\vspace{-14.5pt}

 \begin{gather*}
   \begin{array}{r|cccccccc}
     0\;\,\STRUT \\
     {c_2}\;\,\STRUT & \;\,c_2 \\
     {c_3}\;\,\STRUT & \;\,0 & \;\,c_3 \\
     \frac{1}{2}\;\,\STRUT & \;\,\frac12 - c_2 & \;\,c_2 + c_3 - 1 & 
\;\,1 - c_3 \\
     \frac{1}{2}\;\,\STRUT & \;\,a_{51} & \;\,a_{52} & \;\,a_{53} & 
\phantom{-}a_{54} \\
     {1 - c_3}\;\,\STRUT & \;\,0 & \;\,c_3 & \;\,0 & \phantom{-}a_{64} & 
\phantom{-}a_{65} \\
     {1 - c_2}\;\,\STRUT & \;\,c_2 & \;\,0 & \;\,\frac12 - 2 c_2 & 
\phantom{-}a_{74} & \phantom{-}a_{75} & \;\;\frac12 - 2 c_2 \\
     1\;\,\STRUT & \;\,0 & \;\,c_3 & \;\,0 & \phantom{-}a_{64} & 
\phantom{-}a_{65} & \;\;0 & \phantom{-}c_3 \\
     \hline
     \vphantom{\vert^{\vert^{\vert^\vert}}} & \;\,\frac12 c_2 & 
\;\,\frac12 c_3 & \;\,\frac14 - c_2 & \phantom{-}b_4 & \phantom{-}b_5 & 
\;\;\frac14 - c_2 & \phantom{-}\frac12 c_3 & \phantom{-}\frac12 c_2
   \end{array} \\[5.5pt]
   \left[ \begin{array}{c} \VEC{a}_{4*}\! \\ \VEC{a}_{5*}\! \end{array} 
\right] = \left[ \begin{array}{cccccccc}
     \STRUT{}\PHI c_2 & \phantom{a}(1 - \PHI) c_3 & \phantom{a}\PHI 
(\frac12 - 2 c_2) & \phantom{a}\PHI c_2 + (1 - \PHI) (\frac12 - c_3) & 
\phantom{m}0 & \phantom{-}0 & \phantom{-}0 & \phantom{-}0\, \\
     \STRUT{}\psi c_2 & \phantom{a}(1 - \psi) c_3 & \phantom{a}\psi 
(\frac12 - 2 c_2) & \phantom{a}\psi c_2 + (1 - \psi) (\frac12 - c_3) & 
\phantom{m}0 & \phantom{-}0 & \phantom{-}0 & \phantom{-}0\,
   \end{array} \right] \\[5.5pt] \left[ \begin{array}{cc}
     a_{64} &  \;\,a_{65} \\
     a_{74} & \;\,a_{75} \\
     b_4 & \;\,b_5
   \end{array} \right] = \left[ \begin{array}{cc}
    \STRUT{}2 (\frac12 - c_3) (1 - \chi) & \phantom{-}2 (\frac12 - c_3) 
\chi \\
    \STRUT{}2 c_2 (1 + \chi) & -2 c_2 \chi \\
    \STRUT{}\frac12 (a_{64} + a_{74}) & \phantom{-}\frac12 (a_{65} + 
a_{75})
  \end{array} \right] , \qquad
  \chi = \frac {4 (1 - 3 c_2)} {(1 - 6 c_2) (\PHI - \psi)}
 \end{gather*}

\vspace{-5.5pt}

\NI{}Here $c_2 (c_2 - \frac{1}{2}) (c_2 - 1) = \frac{1}{24}$, $c_3 = c_2 
(1 - c_2) / (\frac12 - c_2) = 1 / 6 (1 - 2 c_2)^2$, and $\PHI = (\frac12 
- c_3) / (\frac12 - c_2 - c_3) = 1 / 2 c_2 - 1$. Out of these three 
points, the methods at the point $(c_2 \SEP c_3) = \mathsf{A} = (z_1 
\SEP z_2) = (0.12888... \SEP 0.30253...)$ are the most efficient.

The residuals in the properties $D \bigl( \VEC{\Phi}(\textrm{t}_1), 
\VEC{\Phi}(\textrm{t}_2) \bigr)$ with $\vert \textrm{t}_1 \vert = 4$ and 
$\vert \textrm{t}_2 \vert = 5$ are inversely proportional to $b_5$. In 
the limit $b_5 \to \infty$, or $\psi \to \PHI$, one gets a method that 
uses information about the derivative of the
 r.h.s.{\ }function $\VEC{f}$ \citep[s.~224]{But16}, 
\citep[s.~IV.7]{HW96}. As the sum $b_4 + b_5$ does not depend on $\psi$, 
maximizing $b_5$ but keeping $b_4 \ge 0$ leads to $b_4 = 0$, $\psi = 
(\frac{1}{2} - c_3) / (\frac{1}{2} + c_2 - c_3) = 2 c_3$:
 \begin{align} \!\!\begin{array}{r|cccccccc}
     0{\;\,}\STRUT & & & & & & \multicolumn{3}{l}{\mbox{\small$c_2 = 
\frac{1}{2} - \frac{1}{\sqrt{3}} \sin \bigl( \mkern-1mu \frac{2 \pi}{9} 
\bigr)$}} \\
     c_2{\;\,}\STRUT & \phantom{l}c_2 & & & & & 
\multicolumn{3}{l}{\mbox{\small$c_3 = \frac{1}{2} - \frac{1}{\sqrt{3}} 
\sin \bigl( \mkern-0.5mu \frac{\pi}{9} \bigr) $}} \\
     c_3{\;\,}\STRUT & \phantom{l}0 & c_3 \\
     \frac{1}{2}{\;\,}\STRUT & \phantom{l}\frac12 - c_2 & c_2 + c_3 - 1 
& 1 - c_3 \\
     \frac{1}{2}{\;\,}\STRUT & \phantom{l}2 c_2 c_3 & (1 - 2 c_3) c_3 & 
(1 - 4 c_2) c_3 & 4 c_2 c_3 \\
     {1 - c_3}{\;\,}\STRUT & \phantom{l}0 & c_3 & 0 & 4 c_2 - 2 & 
\frac{1}{2 c_2} - 2 \\
     {1 - c_2}{\;\,}\STRUT & \phantom{l}c_2 & 0 & \frac12 - 2 c_2 & 2 - 
4 c_2 & 6 c_2 - 2 & \frac12 - 2 c_2 \\
     1{\;\,}\STRUT & \phantom{l}0 & c_3 & 0 & 4 c_2 - 2 & \frac{1}{2 
c_2} - 2 & 0 & \phantom{n}c_3 \\
     \hline
     \vphantom{\vert^{\vert^{\vert^\vert}}} & \phantom{l}\frac12 c_2 & 
\frac12 c_3 & \frac14 - c_2 & 0 & \frac12 + c_2 - c_3 & \frac14 - c_2 & 
\phantom{n}\frac12 c_3 & \phantom{m}\frac12 c_2
   \end{array}\!\! \hskip-12pt \label{the_method} \end{align} All the nodes, 
weights, and coefficients are elements of the algebraic extension 
$\pmb{\mathbf{Q}}(c_2) = \bigl\{ \alpha_1 + \alpha_2 c_2 + \alpha_3 c_3 
\mkern2mu \big\vert \mkern2mu \alpha_1 , \alpha_2 , \alpha_3 \in 
\pmb{\mathbf{Q}} \bigr\}$ of the field of rational numbers 
$\pmb{\mathbf{Q}}$. For example, $c_2^2 = -\frac{1}{12} + \frac{7}{6} 
c_2 - \frac{1}{6} c_3$, $c_2 c_3 = -\frac{1}{12} + \frac{1}{6} c_2 + 
\frac{1}{3} c_3$, $c_3^2 = -\frac{1}{3} + \frac{1}{6} c_2 + \frac{4}{3} 
c_3$, and $a_{65} = a_{85} = 1 / 2 c_2 - 2 = \PHI - 1 = 3 - 4 c_2 - 2 
c_3$.

The stability function of the method in eq.~(\ref{the_method}) is equal 
to

\vspace{-14.5pt}

\begin{gather*}
   R(z) = 1 + z + {\TS\frac{1}{2}} z^2 + {\TS\frac{1}{6}} z^3 + 
{\TS\frac{1}{24}} z^4 + {\TS\frac{1}{120}} r_5 \mkern1mu z^5 + 
{\TS\frac{1}{720}} r_6 \mkern1mu z^6 + {\TS\frac{1}{5040}} r_7 \mkern1mu 
z^7 + {\TS\frac{1}{40320}} r_8 \mkern1mu z^8
 \end{gather*}

\vspace{-3.5pt}

\NI{}where $r_5 = \frac{5}{6} + \frac{10}{3} c_2 - \frac{5}{6} c_3 = 
 1.0108...$, $r_6 = 20 c_2 - 5 c_3 = 6 r_5 - 5 = 1.0650...$, $r_7 = 
 -\frac{35}{8} + 105 c_2 - \frac{105}{4} c_3 = 1.2165...$, and $r_8 = 
 -\frac{70}{3} + \frac{1400}{3} c_2 - \frac{350}{3} c_3 =
1.5179...$

 \begin{table} \vspace{0pt} \begin{tabular}{r|rccclll}
     & $s$ & $p$ & $q$ & $10^4 \times T_4$ & $\phantom{-} 10^3 \times 
T_5$ & $\phantom{-} 10^3 \times T_6$ & $\phantom{-} R(z) R(-z) \,-\, 1$ 
\\
     \hline
     RK4$\STRUT$ & $4$ & $4$ & $4$ & $0$ & $\phantom{0}14.504...$ & 
$\phantom{0}16.035...$ & $\phantom{-}0.01388888... \; z^6$ \\
     AC36 & $5$ & $3$ & $6$ & $7.5690...$ & $\phantom{00}2.3451...$ & 
$\phantom{00}3.9611...$ & $-0.00057870... \; z^8$ \\
     CLMR47$\STRUT$ & $6$ & $4$ & $7$ & $0$ & $\phantom{00}0.29185...$ & 
$\phantom{00}0.39787...$ & $\phantom{-}0.00024760... \; z^8$ \\
     CCRL47 & $6$ & $4$ & $7$ & $0$ & $\phantom{00}1.4813...$ & 
$\phantom{00}1.6278...$ & $-0.00030365... \; z^8$ \\
     eq.~(\ref{method_C})$\STRUT$ & $7$ & $4$ & $9$ & $0$ & 
$112.99...$ & $132.54...$ & $-0.00144678... \; z^{10}$ \\
     eq.~(\ref{the_method}) & $8$ & $4$ & $8$ & $0$ & 
$\phantom{00}0.64048...$ & $\phantom{00}0.91796...$ & $\phantom{-} 
0.00000950... \; z^{10}$ \\
     CV8$\STRUT$ & $11\mkern-0.5mu$ & $8$ & $8$ & $0$ & 
$\phantom{0000}0$ & $\phantom{0000}0$ & $\phantom{-}0.00000627... \; 
z^{10}$ \\
     GL4$\vphantom{\vert_{\vert_\vert}}$ & $2$ & $4$ & $\infty$ & $0$ & 
$\phantom{00}4.3306...$ & $\phantom{00}5.6178...$ & $\phantom{-}0$
 \end{tabular} \begin{tabular}{r|cccccp{4pt}lp{0pt}l}
     & $C(2)\vphantom{\vert^{\big\vert}}$ & 
$\mkern-1mu{}\VEC{D}(\VEC{1})\mkern-1mu$ & $\VEC{D}(\VEC{c})$ & 
$\mkern-3mu{}\VEC{D}(\VEC{c}^2)\mkern-3mu$ & 
$\mkern-7mu{}\VEC{D}(\MAT{A} \mkern0.5mu \VEC{c})\mkern-7mu$ & & 
{$\mbox{max}_{ij} \vert \mkern1mu \AIJ \vert$} & &
{$\phantom{-}\mbox{min}_j \mkern1.5mu \BJ \STRUT$} \\
     \hline
     RK4$\STRUT$ & \FALSE & \TRUE & \FALSE & \FALSE & \FALSE & & 
$\phantom{0}1$ & & $\phantom{-}0.1666...$ \\
     AC36 & \FALSE & \TRUE & \TRUE & \FALSE & \FALSE & & 
$\phantom{0}2.1621...$ & & $-0.3054...$ \\
     CLMR47$\STRUT$ & \FALSE & \TRUE & \TRUE & \TRUE & \TRUE & & 
$\phantom{0}4.4309...$ & & $\phantom{-}0.0277...$ \\
     CCRL47 & \FALSE & \TRUE & \TRUE & \TRUE & \TRUE & & 
$\phantom{0}2.9265...$ & & $-0.4488...$ \\
     eq.~(\ref{method_C})$\STRUT$ & \FALSE & \TRUE & \TRUE & 
\TRUE & \TRUE & & $\phantom{0}1.7024...$ & & $-0.8512...$ \\
     eq.~(\ref{the_method}) & \FALSE & \TRUE & \TRUE & \TRUE & \TRUE & & 
$\phantom{0}1.8793...$ & & $\phantom{-}0.0644...$ \\
     CV8$\STRUT$ & \TRUE & \TRUE & \FALSE & \FALSE & \FALSE & & 
$14.728...$ & & $\phantom{-}0.05$ \\
     GL4$\vphantom{\vert_{\vert_\vert}}$ & \TRUE & \TRUE & \TRUE & \TRUE 
& \TRUE & & $\phantom{0}0.5386... $ & & $\phantom{-}0.5$ \end{tabular} 
\vspace{2pt} \caption{A comparison of eight $s$-stage Runge--Kutta 
methods of order $(p \SEP q)$. Error coefficients are defined as $T_p^2 
= \sum_{\;\textrm{t},\;\vert \mkern1mu \textrm{t} \mkern1mu \vert = p} 
\bigl( \VEC{b} \mkern2mu \VEC{\Phi}(\textrm{t}) - 1 / \textrm{t}! 
\bigr){}^2 \mkern-1mu / \mkern1mu \sigma^2(\textrm{t}) = (1 / p!)^2 
\sum_{\;\textrm{t},\;\vert \mkern1mu \textrm{t} \mkern1mu \vert = p} 
\alpha^2(t) \bigl( \textrm{t}! \mkern3mu \VEC{b} \mkern2mu 
\VEC{\Phi}(\textrm{t}) - 1 \bigr){}^2$, where $\sigma(\textrm{t})$ is 
the order of the symmetry group of the tree $\textrm{t}$, and 
$\alpha(\textrm{t})$ is the number of monotonic labelings of 
$\textrm{t}$ (see, \EG, \cite[ss.~304 and 318]{But16}, \citep[pp.~147 
and 158]{HNW93}, \citep[pp.~58 and 60]{But21}, \citep[pp.~57 and 
58]{HLW06}). The column about the function $R(z) R(-z) - 1$ shows the 
first non-zero term in its Taylor expansion about $z = 0$. Within the 
$1$-dimensional families of methods at the points $\mathsf{A}$, 
$\mathsf{A}'$, and $\mathsf{A}''$ the error coefficients $T_5$, $T_6$, 
and the stability function $R(z)$ do not depend on the parameter $\psi$. 
The symbols \TRUE$\mkern2mu$/$\mkern3mu$\raisebox{0.5pt}{\FALSE} stand 
for true$\mkern1mu$/$\mkern1mu$false. The explicit version of the 
property $C(2)$ is used. The $\mbox{min}_j \mkern1.5mu \BJ$ column shows 
the minimal value of a non-zero weight.} \label{table_methods} 
\end{table}

 \begin{figure} \centerline{\includegraphics{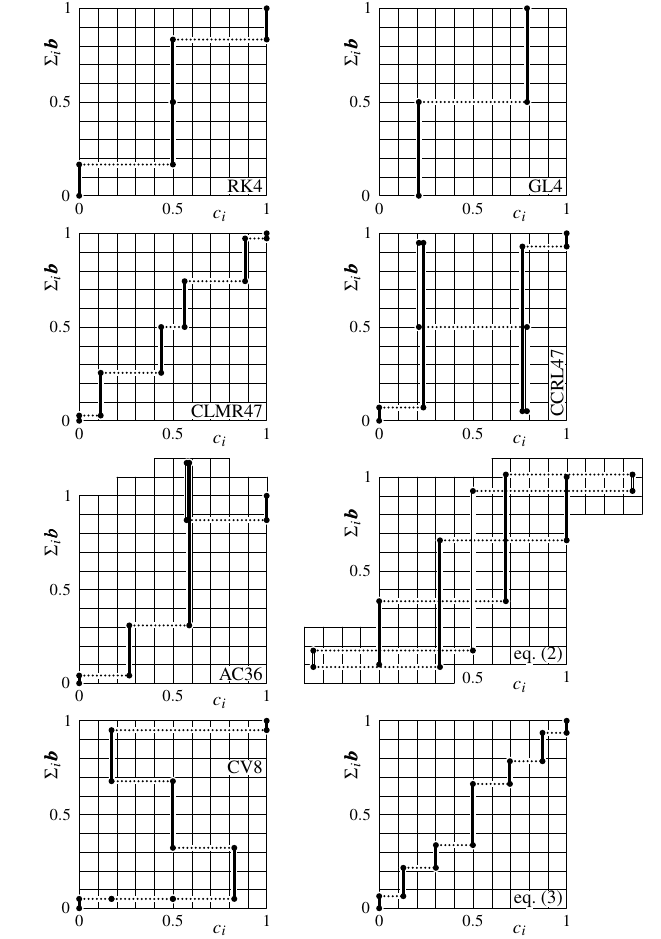}} 
\caption{Quadrature scheme graphical depiction (see also 
\citep[][fig.~1]{Suz90}, where the vertical coordinate is the stage) for 
the eight methods listed in Table~\ref{table_methods}. Here 
${\Sigma_{\mkern1mu i} \mkern1.5mu \VEC{b}} = \sum_{j = 1}^{\mkern1.5mu 
i} \BJ$, where $1 \le i \le s$, are the cumulative weights, with 
$\Sigma_{\mkern1mu 0} \mkern1.5mu \VEC{b} = 0$. Due to the order 
condition $\VEC{b} \mkern0.5mu \VEC{1} = 1$ for an $s$-stage method one 
has $\Sigma_{\mkern1mu s} \mkern1.75mu \VEC{b} = 1$. On the diagram the 
points $\bigl( \CI \SEP {\Sigma_{\mkern1mu i - 1} \mkern1.5mu \VEC{b}} 
\bigr)$ and $\bigl( \CI \SEP {\Sigma_{\mkern1mu i} \mkern1.5mu \VEC{b}} 
\bigr)$, where $1 \le i \le s$, are connected by a thick line (with 
white filling if ${\Sigma_{\mkern1mu i} \mkern1.5mu \VEC{b}} < 
{\Sigma_{\mkern1mu i - 1} \mkern1.5mu \VEC{b}}$). To easier follow the 
progression of the stages, the points $\bigl( \CI \SEP 
{\Sigma_{\mkern1mu i} \mkern1.5mu \VEC{b}} \bigr)$ and $\bigl( 
c_{\mkern0.5mu i + 1} \SEP {\Sigma_{\mkern1mu i} \mkern1.5mu \VEC{b}} 
\bigr)$, where $1 \le i < s$, are connected by dotted lines.} 
\label{stair} \end{figure}

The$\mkern-0.25mu$ basic$\mkern-0.25mu$ properties$\mkern-0.25mu$ 
of$\mkern-0.25mu$ new$\mkern-0.25mu$ methods$\mkern-0.25mu$ 
eq.~(\ref{method_C}), eq.~(\ref{the_method}) and$\mkern-0.25mu$ 
of$\mkern-0.25mu$ some$\mkern-0.25mu$ known$\mkern-0.25mu$ 
Runge--Kutta$\mkern-0.25mu$ methods$\mkern-0.25mu$ are$\mkern-0.25mu$ 
compared$\mkern-0.25mu$ in$\mkern-0.25mu$ Figure~\ref{stair} and in 
Table~\ref{table_methods}, where (and also in Figures~\ref{rigid} 
and$\mkern-0.5mu$ \ref{pendulum}) the$\mkern-0.5mu$ latter$\mkern-0.5mu$ 
methods$\mkern-0.5mu$ are$\mkern-0.5mu$ named$\mkern-0.5mu$ 
as$\mkern-0.5mu$ follows: RK4 is$\mkern-0.5mu$ the$\mkern-0.5mu$ 
classical$\mkern-0.5mu$ Runge--Kutta$\mkern-0.25mu$ 
method$\mkern-0.25mu$ \citep[p.~448]{Kut01}; AC36$\mkern-0.25mu$ 
is$\mkern-0.25mu$ \citep[fig.~4.1]{AuCh98a};$\mkern1mu$\footnote{~In 
\citep[p.~453]{AuCh98a} it is stated that a certain optimization over 
methods satisfying $\hat{D}([\tau])$ (same as property 
$\VEC{D}(\VEC{c})$) results in $c_2 = 0.000050587...$ Nevertheless, 
optimization over methods not necessarily satisfying $\hat{D}([\tau])$ 
led to the method \citep[fig.~4.1]{AuCh98a} that satisfies 
$\hat{D}([\tau])$ and has $c_2 = 0.13502...$} CLMR47 is 
\citep[p.~262]{CLMR10};$\mkern1mu$\footnote{~In \citep[p.~262]{CLMR10} 
the denominators ``355,568'', ``831,328'', and ``9,955,904'' in 
$a_{42}$, $a_{51}$, and $a_{62}$ should be read as ``3,555,680'', 
``8,313,280'', and ``99,559,040'', respectively.} CCRL47 is 
\citep[p.~90]{CCRL17}; CV8 is \citep[tab.~1]{CoVe72}, 
\citep[p.~210]{But16}; and GL4 is the $2$-stage Gauss--Legendre method 
\citep[]{HaHo55}, \citep[p.~56]{But64a}.

\section{Numerical tests} \label{numerical_tests}

\hskip\parindent{}In this section the Runge--Kutta methods from 
Table~\ref{table_methods} are tested on three numerical examples: 
torque-free rotation of a rigid body \citep[s.~36]{LaLi76}, a pendulum 
with a non-separable Hamiltonian \citep[ex.~16.2, p.~313]{HNW93}, and a 
periodic Toda lattice \citep{Tod67}.

\newpage

$\phantom{.}$

\vspace{-20.2pt}

\begin{figure} \centerline{\includegraphics{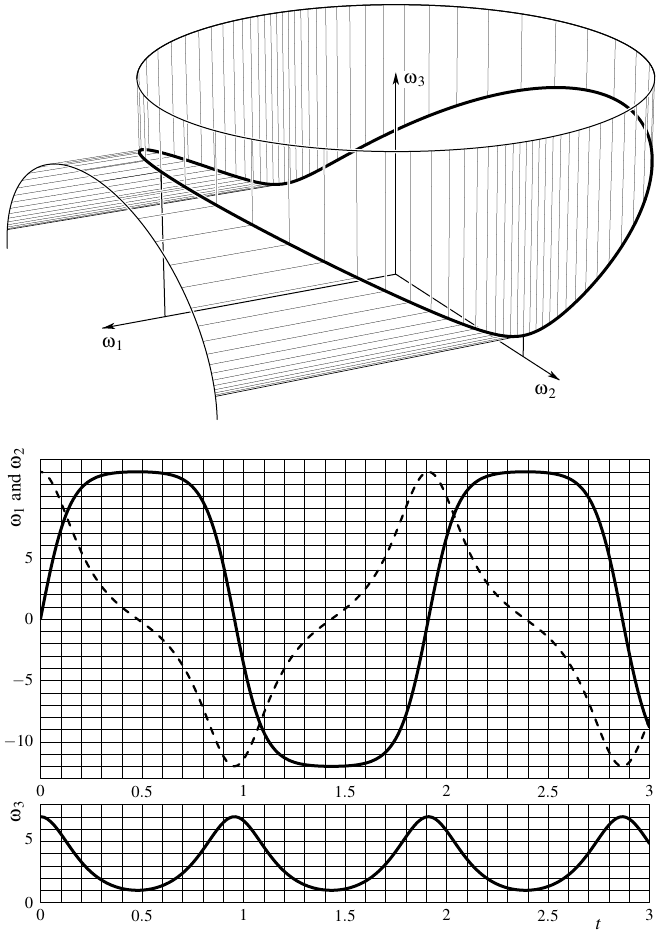}} 
\caption{A solution of Euler's equations with moments of inertia $I_1 = 
1$, $I_2 = 2$, and $I_3 = 3$: $(\OMEGO \SEP \OMEGT \SEP \OMEGH) \BRAT = 
(12 \mkern1mu \textrm{cn} \SEP 12 \mkern1mu \textrm{sn} \SEP 7 \mkern1mu 
\textrm{dn}) \bigl( 7 \mkern0.5mu t \SEP \frac{48}{49} \bigr)$, where 
$\textrm{cn}$, $\textrm{sn}$, and $\textrm{dn}$ are the Jacobi elliptic 
functions. The motion is periodic with the period $4 K \bigl( 
\frac{48}{49} \bigr) / 7 = 1.9109...$, where $K(m) = 
\smash{\int_{\mkern1mu 0}^{\mkern1mu \pi/2}} \d \theta \mkern1.5mu / 
\bigl( 1 - m \mkern1mu \sin^2 \theta \bigr){}^{1/2}$ is the complete 
elliptic integral of the \FIRST{} kind. The circular and elliptic 
cylinders correspond to the quadratic invariants $Q_1(\VEC{\omega}) = 
\OMEGO^{\mkern1mu 2} + \OMEGT^{\mkern1mu 2} = 144$ and 
$Q_2(\VEC{\omega}) = \OMEGT^{\mkern1mu 2} + 3 \OMEGH^{\mkern1mu 2} = 
147$, respectively. The thin lines on them divide the period into $72$ 
equal parts. Points of the trajectory that are connected by vertical 
lines with $\OMEGO$ and $\OMEGT$ axes are $(12 \SEP 0 \SEP 7)$ and $(0 
\SEP 12 \SEP 1)$. The alternating behavior of $\OMEGT\BRAT$ is due to 
the rotation around the intermediate principal axis being unstable 
\cite[p.~48]{Poi34}, \citep[s.~37]{LaLi76}.} \label{csdn} \end{figure}

\begin{figure}
\vspace{-3pt}
   \centerline{\includegraphics{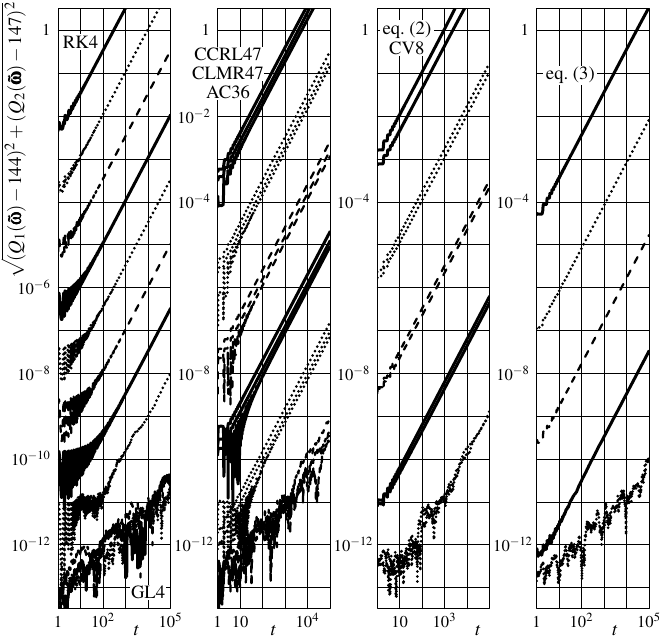}}
\vspace{-3pt}
 \caption{Growth of the error in $Q_1(\VEC{\tilde\omega})$ and 
$Q_2(\VEC{\tilde\omega})$ with time $t$. Panels from left to right: (1) 
RK4 (curves from top to near bottom) and GL4 (bottom curves); (2) AC36 
(lower curves in groups of three), CLMR47 (middle curves), and CCRL47 
(upper curves); (3) CV8 and eq.~(\ref{method_C}) (lower and upper curves 
in pairs); and (4) eq.~(\ref{the_method}). The step size is $h = s 
\mkern1mu h_1$, where $s$ is the number of stages. Upper solid, dotted, 
and dashed curves correspond to $h_1 = 2^{-7}$, $h_1 = 2^{-8}$, and $h_1 
= 2^{-9}$, respectively. Next to upper solid or dotted curves correspond 
to $h_1 = 2^{-10}$ or $h_1 = 2^{-11}$, respectively, and so on, \IE, 
going to the next, with a different line style, curve down divides $h_1$ 
by $2$. Upper solid curves, upper dotted curves, and so on, on different 
panels correspond to the same value of $h_1$, and thus to the same 
number of the r.h.s. function $\VEC{f}(t \SEP \VEC{x})$ evaluations.} 
\label{rigid} \end{figure}

Euler's equations describing free rotation of a rigid body with 
principal moments of inertia $I_1 = 1$, $I_2 = 2$, and $I_3 = 3$ are $\d 
\OMEGO / \DT = -\OMEGT \OMEGH$, $\d \OMEGT / \DT = \OMEGO \OMEGH$, and 
$\d \OMEGH / \DT = -\frac{1}{3} \OMEGO \OMEGT$. The solution with 
initial condition $(\OMEGO \SEP \OMEGT \SEP \OMEGH) \big\vert_{\mkern1mu 
t = 0} = (12 \SEP 0 \SEP 7)$ is shown in Figure~\ref{csdn}. The pair of 
quadratic invariants, twice the kinetic energy $2 T = I_1 
\OMEGO^{\mkern1mu 2} + I_2 \OMEGT^{\mkern1mu 2} + I_3 \OMEGH^{\mkern1mu 
2} = 291$ and squared magnitude of the angular momentum $L^2 = I_1^2 
\OMEGO^{\mkern1mu 2} + I_2^2 \OMEGT^{\mkern1mu 2} + I_3^2 
\OMEGH^{\mkern1mu 2} = 585$, can be reduced to $Q_1(\VEC{\omega}) = 
\OMEGO^{\mkern1mu 2} + \OMEGT^{\mkern1mu 2} = 144$ and 
$Q_2(\VEC{\omega}) = \OMEGT^{\mkern1mu 2} + 3 \OMEGH^{\mkern1mu 2} = 
147$. When this system of differential equations is solved numerically, 
the values of the invariants $Q_1$ and $Q_2$ gradually drift from their 
initial values $144$ and $147$, see Figure~\ref{rigid}. Due to round-off 
errors, $Q_1 \bigl( \VEC{\tilde\omega}(n h) \bigr)$ and $Q_2 \bigl( 
\VEC{\tilde\omega}(n h) \bigr)$ as functions of $n$ can be approximately 
viewed as random walks. In the absence of a systematic drift the typical 
deviation of $(Q_1 \SEP Q_2)$ from $(144 \SEP 147)$ grows roughly 
proportional to $\sqrt{n}$, \EG, the slope of the curves near the bottom 
of Figure~\ref{rigid}, where the round-off errors dominate, is about 
$\frac{1}{2}$. Even a small bias makes the deviation to eventually, \IE, 
for large enough $n$, be proportional to $n$: in the log-log plots of 
Figure~\ref{rigid} the slopes of most of the curves are equal to $1$. 
The performances of AC36, CLMR47, and CCRL47 methods are close, with 
AC36 doing the best out of the three. Besides the implicit GL4, the 
methods in eq.~(\ref{method_C}) and eq.~(\ref{the_method}), and also 
CV8, due to their higher pseudo-symplecticity order have the fastest 
rate of decrease of the error in $Q_1$ and $Q_2$ with the step size $h$ 
(see the two right panels in Figure~\ref{rigid}). The method in 
eq.~(\ref{the_method}), out of the explicit methods tested, has the best 
performance.

 \begin{figure}
   \vspace{10pt}
   \centerline{\includegraphics{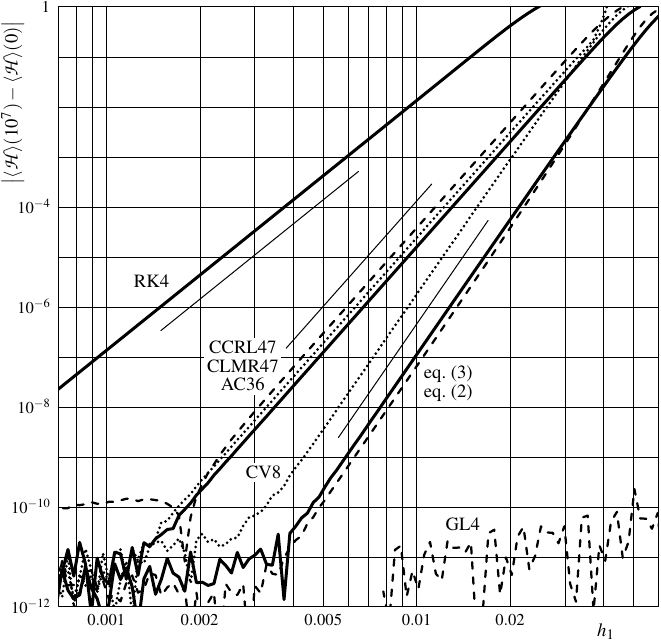}}
 \caption{The speed of the systematic drift of the Hamiltonian 
$\mathcal{H}$ as the function of the step size $h = s \mkern1mu h_1$, 
where $s$ is the number of stages. The number of the r.h.s. function 
$\VEC{f}(t \SEP \VEC{x})$ evaluations is equal to the time duration of 
the simulation (which is $10020000$ here) divided by $h_1$, and with 
$h_1$ being fixed is independent of the method. The curves correspond to 
RK4, AC36, and eq.~(\ref{the_method}) (solid curves from top to bottom), 
CLMR47 and CV8 (upper and lower dotted curves), and CCRL47, 
eq.~(\ref{method_C}), and GL4 (dashed curves from top to bottom). Three 
thin solid lines on this log-log plot have slopes $5$, $7$, and $9$, 
which corresponds to the speed of the drift being proportional to $h^5$, 
$h^7$, and $h^9$. The weighted moving average of the Hamiltonian that 
was used is $\langle \mathcal{H} \rangle \BRAT = \sum_{\mkern1.5mu n 
\mkern1.5mu \in \mathcal{S} \BRAT} \sin^2 \bigl( \pi (n h - t) / 20000 
\bigr) \mkern1mu \mathcal{H} \bigl( \tilde{p}(n h) \SEP \tilde{x}(n h) 
\bigr) \mkern1mu \big/ \mkern1mu \sum_{\mkern1.5mu n \mkern1.5mu \in 
\mathcal{S} \BRAT} \sin^2 \bigl( \pi (n h - t) / 20000 \bigr)$, where 
$\mathcal{S} \BRAT = \bigl\{ n \mkern2.5mu \big\vert \mkern3.5mu t < n h 
< t + 20000 \bigr\}$. The simulations were run at the University of 
Arizona High Performance Computing center.} \label{pendulum} 
\end{figure}

$\phantom{.}$

\vspace{-15.5pt}

\begin{figure}
   \vspace{10pt}
   \centerline{\includegraphics{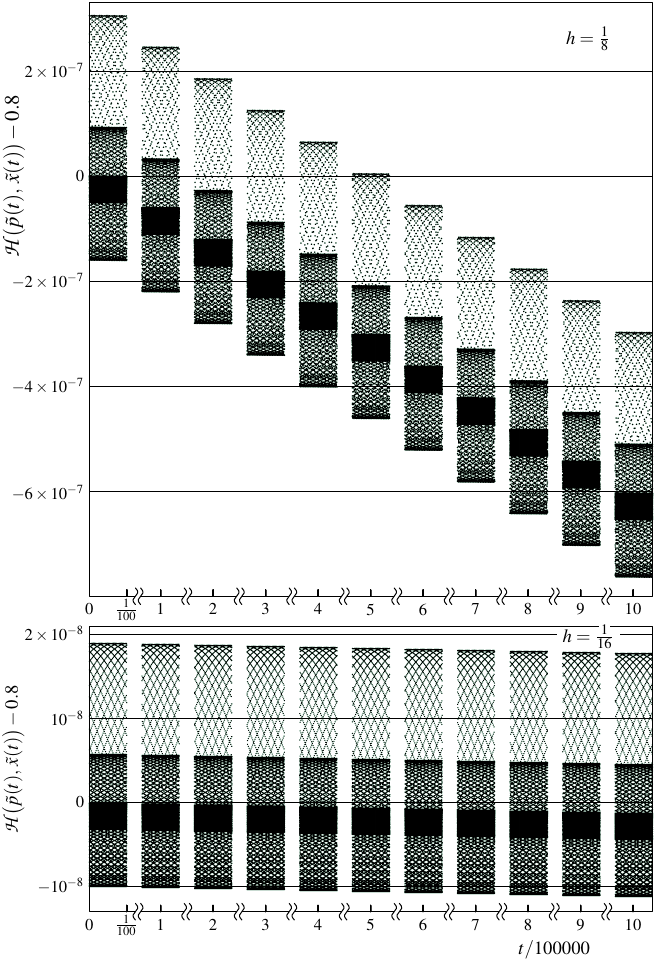}}
 \caption{The Hamiltonian $\mathcal{H} \bigl( \tilde{p}(t) \SEP 
\tilde{x}(t) \bigr)$ computed on the numerical solution obtained by the 
method in eq.~(\ref{the_method}). The vertical scale in the lower panel 
(time step $h = \frac{1}{16}$) is $16$ times larger than in the upper 
one ($h = \frac{1}{8}$). Similarity of heights of the depicted 
Hamiltonian deviations in the two panels is the consequence of the 
method being of order $4$. In order to show both the deviations small 
time scale structure and the systematic drift of the Hamiltonian, from 
the whole simulation (which lasts $1000500$) eleven fragments of 
duration $1000$ are shown. As the method conserves the symplectic 
structure up to the order $8$, the speed of the systematic drift 
decreases with $h$ much faster than the magnitude of the deviations 
within a period.} \label{linear} \end{figure}

The \SECOND{} numerical example is a system $\d \mkern0.5mu x / \DT = 
\partial \mathcal{H} / \partial p$, $\d p / \DT = -\partial \mathcal{H} 
/ \partial x$, with the Hamiltonian function $\mathcal{H}(p \SEP x) = 
\frac{1}{2} p^2 - \bigl( 1 - \frac{1}{6} p \bigr) \cos(x)$. The solution 
with initial condition $(x \SEP p) \big\vert_{\mkern1mu t = 0} = \bigl( 
\arccos(-0.8) \SEP 0 \bigr)$ forms a periodic trajectory along the 
closed level curve $\mathcal{H}(p \SEP x) = 0.8$ around the origin. 
Because of notable variations of the Hamiltonian within a period in a 
numerical solution (see, \EG, Figure~\ref{linear}), here the speed of 
the systematic drift of the Hamiltonian $\mathcal{H}$ is estimated 
through the difference between weighted moving average values at two 
different moments of time. For the eight methods in 
Table~\ref{table_methods}, how the drift of the Hamiltonian depends on 
the step size $h$, is shown in Figure~\ref{pendulum}. Again, the 
performances of AC36, CLMR47, and CCRL47 methods are close, with AC36 
doing the best out of the three. In$\mkern-0.2mu$ the$\mkern-0.2mu$ 
triple$\mkern-0.2mu$ of$\mkern-0.2mu$ higher$\mkern-0.2mu$ 
order$\mkern-0.2mu$ explicit$\mkern-0.2mu$ methods, the$\mkern-0.2mu$ 
ones$\mkern-0.2mu$ in$\mkern-0.2mu$ eq.~(\ref{method_C})$\mkern-0.2mu$ 
and$\mkern-0.2mu$ eq.~(\ref{the_method})$\mkern-0.2mu$ 
outperform$\mkern-0.2mu$ CV8.\footnote{~Although$\mkern-0.8mu$ 
with$\mkern-0.8mu$ the$\mkern-0.8mu$ same$\mkern-0.8mu$ 
amount$\mkern-0.8mu$ of$\mkern-0.8mu$ work$\mkern-0.8mu$ 
the$\mkern-0.8mu$ drift$\mkern-0.8mu$ of$\mkern-0.8mu$ the$\mkern-0.8mu$ 
Hamiltonian$\mkern-0.8mu$ for$\mkern-0.8mu$ CV8$\mkern-0.8mu$ 
is$\mkern-0.8mu$ about$\mkern-0.8mu$ $15${}$\mkern-0.8mu$ 
times$\mkern-0.8mu$ faster$\mkern-0.8mu$ than$\mkern-0.8mu$ 
for$\mkern-0.8mu$ the$\mkern-0.8mu$ method$\mkern-0.8mu$ 
in$\mkern-0.8mu$ eq.~(\ref{the_method}), the$\mkern-0.8mu$ 
course$\mkern-0.8mu$ of$\mkern-0.8mu$ trajectory$\mkern-0.8mu$ $\bigl( 
x(t) \SEP p(t) \bigr)${}$\mkern-0.8mu$ within$\mkern-0.8mu$ 
a$\mkern-0.8mu$ period$\mkern-0.8mu$ is$\mkern-0.8mu$ 
computed$\mkern-0.8mu$ more$\mkern-0.8mu$ accurately$\mkern-0.8mu$ 
by$\mkern-0.8mu$ the$\mkern-0.8mu$ \NTH{8}$\mkern-0.8mu$ 
order$\mkern-0.8mu$ method$\mkern-0.8mu$ CV8.}

The \THIRD{} numerical example is the periodic Toda lattice with $N = 
32$ nodes:

\vspace{-17pt}

 \begin{gather*}
   \mathcal{H}(\VEC{p} \SEP \VEC{x}) = \frac{1}{2} \sum_{n = 1}^N p_n^2 
+ \sum_{n = 1}^N U \bigl( x_n - x_{n - 1} \bigr) , \qquad U(r) = 
\exp(-r) + r - 1 \\
   \frac{\d \mkern0.5mu x_n}{\DT} = \frac{\partial \mathcal{H}}{\partial 
p_n} = p_n , \qquad
   \frac{\d^2 x_n}{\DT^2} = \frac{\d p_n}{\DT} = -\frac{\partial 
\mathcal{H}}{\partial x_n} = \exp(x_{n - 1} - x_n) - \exp(x_n - x_{n + 
1}) \end{gather*}

\vspace{-6pt}

\NI{}with the periodicity condition $x_0 \equiv x_N$. This Hamiltonian 
system is completely integrable and can be reformulated as the Lax 
equation $\d \mkern1.5mu \MAT{L} / \DT = [ \mkern1.5mu \MAT{L} \SEP 
\MAT{A} ] = \MAT{L} \mkern2mu \MAT{A} - \MAT{A} \mkern1mu \MAT{L}$, with 
the Lax pair
 \begin{gather*}
   \\[-16pt]
   \MAT{L} = \left[ \begin{array}{ccccccc}
     b_1 & a_1 & 0 & 0 & \cdots & 0 & a_N \\
     \,a_1\, & b_2 & \,a_2\, & 0 & \cdots & 0 & 0 \\
     0 & \,a_2\, & b_3 & \,a_3\, & \cdots & 0 & 0 \\
     0 & 0 & a_3 & b_4 & \cdots & 0 & 0 \\
     \vdots & \vdots & \vdots & \vdots & \ddots & \vdots & \vdots \\
     0 & 0 & 0 & 0 & \cdots & \!b_{N - 1}\! & \!a_{N - 1}\! \\
     a_N & 0 & 0 & 0 & \cdots & \!a_{N - 1}\! & \!b_N\! \end{array} \right] , \quad\!\!
   \MAT{A} = \left[ \begin{array}{ccccccc}
     0 & \!\!\!\!\!-a_1 & 0 & 0 & \cdots & 0 & a_N \\
     \,\,a_1\,\, & 0 & \!\!\!\!\!-a_2 & 0 & \cdots & 0 & 0 \\
     0 & \,\,a_2\,\, & 0 & \!\!\!-a_3\,\, & \cdots & 0 & 0 \\
     0 & 0 & \,\,a_3\,\, & 0 & \cdots & 0 & 0 \\
     \vdots & \vdots & \vdots & \vdots & \ddots & \vdots & \vdots \\
     0 & 0 & 0 & 0 & \cdots & 0 & \!\!\!\!\!\!-a_{N - 1}\! \\
     \!\!\!\!-a_N\, & 0 & 0 & 0 & \cdots & a_{N - 1} & 0 \end{array} \right]
   \\[-16pt]
 \end{gather*} where $a_n = {\TS\frac12} \exp \bigl( {\TS\frac12} (x_{n 
- 1} - x_n) \bigr)$, $b_n = -{\TS\frac12} p_{n - 1}$ are the so-called 
Flaschka variables \citep{Fla74}. Due to periodicity, $p_0 \equiv p_N$ 
and $b_1 = -{\TS\frac12} p_N$. All $N$ eigenvalues $\lambda_1$, 
$\lambda_2$, \DOTS, $\lambda_N$ of $\MAT{L}$ are integrals of motion 
\citep{Lax68}. The coefficients of the characteristic polynomial $\det 
(\lambda \MAT{I} - \MAT{L})$ are related to the integrals of motion 
found in \citep{Hen74}, see \citep{Fla74} for details. In particular, the 
total momentum $P = \sum_{n = 1}^N p_n = -2 \mkern1mu \textrm{tr} \, 
\MAT{L} = -2 \sum_{n = 1}^N \lambda_n$ and the Hamiltonian $\mathcal{H} 
= 2 \mkern1mu \textrm{tr} \, \MAT{L}^2 - N = 2 \sum_{n = 1}^N 
\lambda_n^2 - N$ do not depend on time.

\begin{figure}
   \centerline{\includegraphics{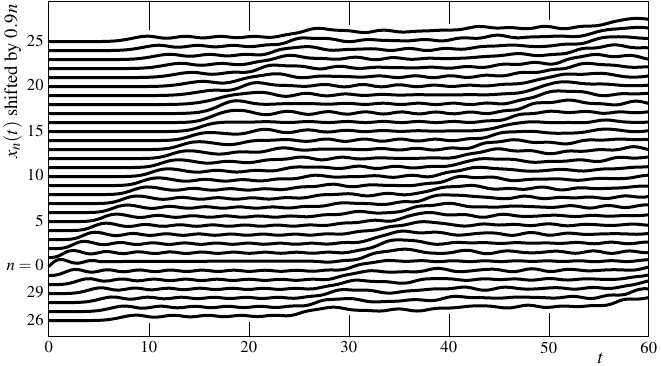}}
 \caption{Evolution of the periodic Toda lattice with $N = 32$ and 
initial condition $x_n(0) = 0$ for all $1 \le n \le 32$, $p_n(0) = 0$ 
for all $1 \le n < 32$, and $p_0(0) = p_{32}(0) = 1$. To show wave-like 
perturbations propagating both to the right and to the left from the 
initial disturbance at $n = 32 \equiv 0$, the curves corresponding to $n 
\ge 26$ are drawn at the bottom, with $n$ being treated modulo $32$. All 
the curves slowly but systematically climb up, as the total momentum $P 
= \sum_{n = 1}^N p_n(t) = 1$ is positive.} \label{evolution} 
\end{figure}

The initial condition $x_n(0) = 0$ for all $1 \le n \le N$, $p_n(0) = 0$ 
for all $1 \le n < N$, and $p_N(0) = 1$ was used in the simulations. It 
corresponds to $P = 1$, $\mathcal{H} = {\TS\frac12}$, $a_n(0) = 
{\TS\frac12}$ for all $1 \le n \le N$, $b_n(0) = 0$ for all $1 \le n < 
N$, and $b_N(0) = -{\TS\frac12}$. The whole system is in equilibrium and 
at rest, with the exception of the \NTH{N} particle which moves with 
velocity $p_N(0) = 1$. The evolution of the system over the time 
interval $0 \le t \le 60$ is shown in Figure~\ref{evolution}.

\begin{figure}[h]
   \centerline{\includegraphics{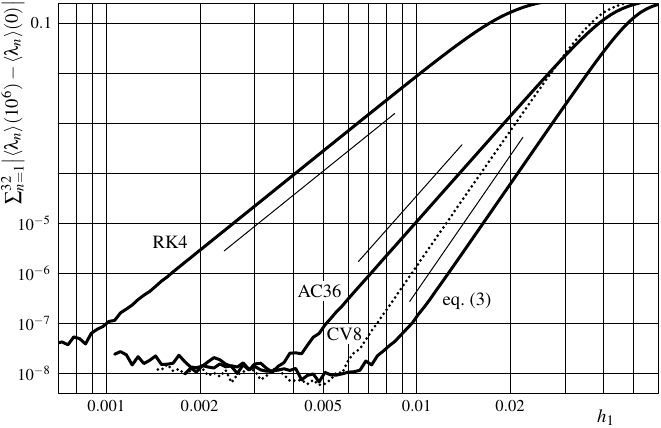}}
 \caption{The combibed speed of the systematic drift of $N = 32$ 
integrals of motion $\lambda_n$, $1 \le n \le 32$, that are the 
eigenvalues of $\MAT{L}$ as the function of the step size $h = s 
\mkern1mu h_1$, where $s$ is the number of stages. The number of the 
r.h.s. function $\VEC{f}(t \SEP \VEC{x})$ evaluations is equal to the 
time duration of the simulation (which is $1002000$ here) divided by 
$h_1$, and with $h_1$ being fixed is independent of the method. The 
curves correspond to RK4, AC36, and eq.~(\ref{the_method}) (solid curves 
from top to bottom), and CV8 (dotted curve). Three thin solid lines on 
this log-log plot have slopes $5$, $7$, and $9$, which corresponds to 
the speed of the drift being proportional to $h^5$, $h^7$, and $h^9$. 
The weighted moving average that was used is $\langle \lambda \rangle 
\BRAT = \sum_{\mkern1.5mu n \mkern1.5mu \in \mathcal{S} \BRAT} \sin^2 
\bigl( \pi (n h - t) / 2000 \bigr) \mkern1mu \lambda(n h) \mkern1mu 
\big/ \mkern1mu \sum_{\mkern1.5mu n \mkern1.5mu \in \mathcal{S} \BRAT} 
\sin^2 \bigl( \pi (n h - t) / 2000 \bigr)$, where $\mathcal{S} \BRAT = 
\bigl\{ n \mkern2.5mu \big\vert \mkern3.5mu t < n h < t + 2000 \bigr\}$. 
The simulations were run at the University of Arizona High Performance 
Computing center.} \label{lattice} \end{figure}

How the systematic drift of the integrals of motion $\lambda_n$, $1 \le 
n \le 32$, depends on the step size $h$, is shown in 
Figure~\ref{lattice}. The performances of AC36, CLMR47, and CCRL47 
methods are close, with AC36 doing the best out of the three, the CLMR47 
and CCRL47 curves are omitted. The CLMR47 and eq.~(\ref{method_C}) 
methods are prone to cause numerical overflow at large time steps $h$, 
most probably due to the presence of exponential function in the 
r{.}h{.}s. The method in eq.~(\ref{the_method}) is capable of conserving 
up to the \NTH{8} order multiple integrals of motion that are far from 
being quadratic, with a bit better performance than CV8 for the same 
amount of r.h.s. function evaluations.

\section*{Conclusions}

\hskip\parindent{}There are known $5$- and $6$-stage pseudo-symplectic 
Runge--Kutta methods of order $(3 \SEP 6)$ and $(4 \SEP 7)$, 
respectively, see, \EG, Table~\ref{table_methods}. With $7$ stages it is 
possible to come up with a method of order $(4 \SEP 9)$, see 
eq.~(\ref{method_C}). Utilising $8$ stages, one can construct robust, 
with non-negative weights and monotonically increasing nodes, methods of 
order $(4 \SEP 8)$, see eq.~(\ref{the_method}). The newly derived 
methods, largely due to their higher order, have better quadratic 
invariants and energy preservation properties than previously known 
pseudo-symplectic methods.

\newpage

\newpage

$\phantom{.}$

\vspace{-22pt}

\section*{Appendix}

Below is Wolfram Mathematica script that validates the three families of 
Runge--Kutta methods corresponding to the points $(c_2 \SEP c_3) = 
\mathsf{A}$, $\mathsf{A}'$, and $\mathsf{A}''$, and indexed by a 
parameter $\psi$. The output is zero matrix and zero vector to confirm 
the expressions, four zero vectors standing for $\VEC{D}{\BRA1}$, 
$\VEC{D}(\VEC{c})$, $\VEC{D}(\VEC{c}^2)$, and $\VEC{D}(\MAT{A} 
\mkern0.5mu \VEC{c})$, plus vector $\bigl[ \; 1 ~~~ \frac{1}{2} ~~~ 
\frac{1}{3} ~~~ \frac{1}{6} ~~~ \frac{1}{4} ~~~ \frac{1}{8} ~~ 
\frac{1}{12} ~~~ \frac{1}{24} \; \bigr]$, to show that the order 
conditions are satisfied. The script exactly follows the derivation 
steps in Section~\ref{family}.

\bigskip

{\footnotesize\begin{verbatim}
A = {{  0,   0,   0,   0,   0,   0,   0,   0},
     { c2,   0,   0,   0,   0,   0,   0,   0},
     {a31, a32,   0,   0,   0,   0,   0,   0},
     {a41, a42, a43,   0,   0,   0,   0,   0},
     {a51, a52, a53, a54,   0,   0,   0,   0},
     {a61, a32,   0, a64, a65,   0,   0,   0},
     {a71,   0, a73, a74, a75, a73,   0,   0},
     {a81, a82, a83, a84, a85, a86, a87,   0}}
b =  { b1,  b2,  b3,  b4,  b5,  b3,  b2,  b1}
c =  {  0,  c2,  c3, 1/2, 1/2, 1-c3, 1-c2, 1}
u =  {  1,   1,   1,   1,   1,   1,   1,   1}
SOLVE[e_, v_] := Set @@@ Simplify[Solve[e, v]][[1]]
Ap1 = Simplify[A.(2*c - u)]
SOLVE[Ap1[[2]] == Ap1[[7]], a71];  SOLVE[Ap1[[3]] == Ap1[[6]], a61]
Dp1 = Simplify[(b*(2*c - u)).A + b*Ap1]
SOLVE[Dp1[[1]] == 0, a81];  SOLVE[Dp1[[2]] == 0, a82]
SOLVE[Dp1[[7]] == 0, a87];  SOLVE[Ap1[[8]] == 0, a86]
SOLVE[(A.u)[[3]] == c[[3]], a31]
D1 = Simplify[b.A + b*(c - u)]
SOLVE[{D1[[6]] == 0, Dp1[[6]] == 0}, {a73, a83}]
SOLVE[(A.u)[[6]] == c[[6]], a64]
Ap2 = Simplify[A.(6*c^2 - 6*c + u)];  SOLVE[Ap2[[3]] + Ap2[[6]] == 0, a32]
SOLVE[(A.u)[[4]] == c[[4]], a41];  SOLVE[Ap2[[4]] == 0, a42]
Aq1 = Simplify[A.(A.c - (c^2)/2)];  SOLVE[Aq1[[4]] == 0, a43]
SOLVE[(A.u)[[5]] == c[[5]], a51];  SOLVE[Ap2[[5]] == 0, a54]
SOLVE[(A.u)[[7]] == c[[7]], a74];  SOLVE[Ap2[[2]] + Ap2[[7]] == 0, b3]
SOLVE[(A.u)[[8]] == c[[8]], a84];  SOLVE[Ap2[[8]] == 0, b1]
SOLVE[D1[[5]] == 0, a85];  SOLVE[b.u == 1, b4];  SOLVE[b.(c^2) == 1/3, b2]
SOLVE[Dp1[[5]] == 0, a75];  SOLVE[D1[[2]] == 0, a52];  SOLVE[D1[[3]] == 0, a53]
Dc2 = Simplify[(b*(c^2)).A + b*(A.(c^2) - u/3)];  SOLVE[Dc2[[5]] == 0, a65]
SOLVE[A[[5, 1]] == psi*c2, b5]
chi = 4*(1 - 3*c2) / ((1 - 6*c2)*(1/(2*c2) - 1 - psi))
z64 = 2*(1/2 - c3)*(1 - chi);  z65 = 2*(1/2 - c3)*chi
z74 = 2*c2*(1 + chi);  z75 = -2*c2*chi
Z = {{0, 0, 0, 0, 0, 0, 0, 0}, {c2, 0, 0, 0, 0, 0, 0, 0},
  {0, c3, 0, 0, 0, 0, 0, 0}, {1/2 - c2, c2 + c3 - 1, 1 - c3, 0, 0, 0, 0, 0},
  {psi*c2, (1 - psi)*c3, psi*(1/2 - 2*c2), psi*c2 + (1 - psi)*(1/2 - c3),
  0, 0, 0, 0}, {0, c3, 0, z64, z65, 0, 0, 0}, {c2, 0, 1/2 - 2*c2, z74, z75,
  1/2 - 2*c2, 0, 0}, {0, c3, 0, z64, z65, 0, c3, 0}}
y = {c2, c3, 1/2 - 2*c2, z64 + z74, z65 + z75, 1/2 - 2*c2, c3, c2}/2
Print["{A - Z, b - y, D1, Dp1, Dc2, DAc} = ",
  Simplify[{A - Z, b - y, D1, Dp1, Dc2, (b*(A.c)).A + b*(A.(A.c) - u/6)},
  Assumptions->{c2*(c2 - 1/2)*(c2 - 1) == 1/24, 6*c3*(1 - 2*c2)^2 == 1}]]
ORD = {u, c, c^2, A.c, c^3, c*(A.c), A.(c^2), A.(A.c)}
Print["ORD = ", Simplify[b.Transpose[ORD]]]
\end{verbatim}}

\medskip\smallskip

\NI{}The expressions for $\MAT{A}$ and $\VEC{b}$ are confirmed in such a 
fashion because, \EG, even assuming that $c_2 (c_2 - \frac{1}{2}) (c_2 - 
1) = \frac{1}{24}$ and $6 c_3 (1 - 2 c_2)^2 = 1$, due to details of 
implementation, the coefficient $a_{32} = (6 c_3 - 1) / 24 c_2 (1 - 
c_2)$ is unchanged by $\mkern2mu${\tt\bfseries Simplify}$\mkern2mu$ and 
$\mkern2mu${\tt\bfseries FullSimplify}$\mkern2mu$ commands, while 
$a_{32} - c_3$ is correctly simplified to $0$.

\end{document}